\pageno=1
\def\title{Projective completions of Jordan pairs} 
\def\date{17.6.2003}


\input amssym.def
\input amssym.tex

\def\item#1{\vskip1.3pt\hang\textindent {\rm #1}}


\newskip\litemindent
\litemindent=0.7cm  
\def\Litem#1#2{\par\noindent\hangindent#1\litemindent
\hbox to #1\litemindent{\hfill\hbox to \litemindent
{\ninerm #2 \hfill}}\ignorespaces}
\def\litem{\Litem1}

\tolerance=300
\pretolerance=200
\hfuzz=1pt
\vfuzz=1pt

\hoffset=0in
\voffset=0.5in

\hsize=5.8 true in 
\vsize=9.2 true in
\parindent=25pt
\mathsurround=1pt
\parskip=1pt plus .25pt minus .25pt
\normallineskiplimit=.99pt

\countdef\revised=100
\mathchardef\emptyset="001F 
\chardef\ss="19
\def\3{\ss}
\def\anf{$\lower1.2ex\hbox{"}$}
\def\frac#1#2{{#1 \over #2}}
\def\>{>\!\!>}
\def\<{<\!\!<}

\def\into{\hookrightarrow}
 
\def\ssssarr{\hbox to 15pt{\rightarrowfill}}
\def\sssarr{\hbox to 20pt{\rightarrowfill}}
\def\ssarr{\hbox to 30pt{\rightarrowfill}}
\def\sarr{\hbox to 40pt{\rightarrowfill}}
\def\arr{\hbox to 60pt{\rightarrowfill}}
\def\larr{\hbox to 60pt{\leftarrowfill}}
\def\Arr{\hbox to 80pt{\rightarrowfill}}

\def\ad{\mathop{\rm ad}\nolimits}

\def\Aut{\mathop{\rm Aut}\nolimits}

\def\der{\mathop{\rm der}\nolimits}

\def\End{\mathop{\rm End}\nolimits}

\def\GL{\mathop{\rm GL}\nolimits}
\def\Gr{\mathop{\rm Gr}\nolimits}

\def\Hom{\mathop{\rm Hom}\nolimits}%
\def\id{\mathop{\rm id}\nolimits} 
\def\im{\mathop{\rm im}\nolimits}


\def\Pol{\mathop{\rm Pol}\nolimits}




\def\trile{\trianglelefteq}

\def\0{{\bf 0}}
\def\1{{\bf 1}}

\def\e{{\frak e}}
\def\f{{\frak f}}
\def\g{{\frak g}}
\def\gl{{\frak {gl}}}
\def\h{{\frak h}}

\def\q{{\frak q}}

\def\z{{\frak z}}

\def\C{{{\Bbb C}{\mskip+2mu}}} 
\def\K{{{\Bbb K}{\mskip+2mu}}} 

\def\R{{\Bbb R}} 
 
\def\N{{\Bbb N}}

\def\K{{\Bbb K}}

\def\P{{\Bbb P}} 
\def\Q{{\Bbb Q}}

\def\:{\colon}  
\def\.{{\cdot}}
\def\|{\Vert}
\def\bsk{\bigskip}

\def\giantskip{\vskip2\bigskipamount}
\def\gsk{\giantskip}
\def \la {\langle}

\def\msk{\medskip}
\def \ra {\rangle}
\def \res {\!\mid\!\!}

\def\ssk{\smallskip}

\def\bbr{\bigbreak}
\def\giantbreak{\par \ifdim\lastskip<2\bigskipamount \removelastskip
         \penalty-400 \giantskip\fi}

\def\nin{\noindent}
\def\cen{\centerline}
\def\pagebreak{\vskip 0pt plus 0.0001fil\break}
\def\linebreak{\break}

\def\hat{\widehat}

\def\epsilon{\varepsilon}

\def\nin{\noindent}

\def\uline{\underline}
\def\pder#1,#2,#3 { {\partial #1 \over \partial #2}(#3)}
\def\pde#1,#2 { {\partial #1 \over \partial #2}}
\def\phi{\varphi}


\def\subeq{\subseteq}

\def\Rarrow{\Rightarrow}

\def\tilde{\widetilde}

\font\ninerm=cmr9
\font\eightrm=cmr8

\font\eightbf=cmbx8


\font\smc=cmcsc10
\font\bfone=cmbx10 scaled\magstep1 
\font\bftwo=cmbx10 scaled\magstep2 

\def\qed{{\unskip\nobreak\hfil\penalty50\hskip .001pt \hbox{}\nobreak\hfil
          \vrule height 1.2ex width 1.1ex depth -.1ex
           \parfillskip=0pt\finalhyphendemerits=0\medbreak}\rm}

\def\qeddis{\eqno{\vrule height 1.2ex width 1.1ex depth -.1ex} $$
                   \medbreak\rm}

\def\Lemma #1. {\bigbreak\vskip-\parskip\noindent{\bf Lemma #1.}\quad\it}

\def\Sublemma #1. {\bigbreak\vskip-\parskip\noindent{\bf Sublemma #1.}\quad\it}

\def\Proposition #1. {\bigbreak\vskip-\parskip\noindent{\bf Proposition #1.}
\quad\it}

\def\Corollary #1. {\bigbreak\vskip-\parskip\nin{\bf Corollary #1.}
\quad\it}

\def\Theorem #1. {\bigbreak\vskip-\parskip\noindent{\bf Theorem #1.}
\quad\it}

\def\Definition #1. {\rm\bigbreak\vskip-\parskip\noindent{\bf Definition #1.}
\quad}

\def\Remark #1. {\rm\bigbreak\vskip-\parskip\noindent{\bf Remark #1.}\quad}

\def\Example #1. {\rm\bigbreak\vskip-\parskip\noindent{\bf Example #1.}\quad}
\def\Examples #1. {\rm\bigbreak\vskip-\parskip\noindent{\bf Examples #1.}\quad}

\def\Problems #1. {\bigbreak\vskip-\parskip\noindent{\bf Problems #1.}\quad}
\def\Problem #1. {\bigbreak\vskip-\parskip\noindent{\bf Problem #1.}\quad}
\def\Exercise #1. {\bigbreak\vskip-\parskip\noindent{\bf Exercise #1.}\quad}

\def\Conjecture #1. {\bigbreak\vskip-\parskip\noindent{\bf Conjecture #1.}\quad}

\def\Proof#1.{\rm\par\ifdim\lastskip<\bigskipamount\removelastskip\fi\smallskip
            \noindent {\bf Proof.}\quad}

\def\Axiom #1. {\bigbreak\vskip-\parskip\noindent{\bf Axiom #1.}\quad\it}

\def\Satz #1. {\bigbreak\vskip-\parskip\noindent{\bf Satz #1.}\quad\it}

\def\Korollar #1. {\bbr\vskip-\parskip\nin{\bf Korollar #1.} \quad\it}

\def\Folgerung #1. {\bbr\vskip-\parskip\nin{\bf Folgerung #1.} \quad\it}

\def\Folgerungen #1. {\bbr\vskip-\parskip\nin{\bf Folgerungen #1.} \quad\it}

\def\Bemerkung #1. {\rm\bigbreak\vskip-\parskip\noindent{\bf Bemerkung #1.}
\quad}

\def\Beispiel #1. {\rm\bigbreak\vskip-\parskip\noindent{\bf Beispiel #1.}\quad}
\def\Beispiele #1. {\rm\bigbreak\vskip-\parskip\noindent{\bf Beispiele #1.}\quad}
\def\Aufgabe #1. {\rm\bigbreak\vskip-\parskip\noindent{\bf Aufgabe #1.}\quad}
\def\Aufgaben #1. {\rm\bigbreak\vskip-\parskip\noindent{\bf Aufgabe #1.}\quad}

\def\Beweis#1. {\rm\par\ifdim\lastskip<\bigskipamount\removelastskip\fi
           \smallskip\noindent {\bf Beweis.}\quad}

\nopagenumbers

\def\date{\ifcase\month\or January\or February \or March\or April\or May
\or June\or July\or August\or September\or October\or November
\or December\fi\space\number\day, \number\year}

\def\title{Title ??}
\def\author{Author ??}

\def\thanks#1{\footnote*{\eightrm#1}}

\def\rightheadline{\hfil{\eightrm\title}\hfil\tenbf\folio}
\def\leftheadline{\tenbf\folio\hfil{\eightrm\author}\hfil}
\headline={\vbox{\line{\ifodd\pageno\rightheadline\else\leftheadline\fi}}}

\def\firstheadline{}
\def\firstfootline{\cen{\rm\folio}}

\def\seite #1 {\pageno #1
               \headline={\ifnum\pageno=#1 \firstheadline
               \else\ifodd\pageno\rightheadline\else\leftheadline\fi\fi}
               \footline={\ifnum\pageno=#1 \firstfootline\else{}\fi}}

\newdimen\dimenone
 \def\checkleftspace#1#2#3#4{
 \dimenone=\pagetotal
 \advance\dimenone by -\pageshrink   
 \ifdim\dimenone>\pagegoal          
   \else\dimenone=\pagetotal
        \advance\dimenone by \pagestretch
        \ifdim\dimenone<\pagegoal
          \dimenone=\pagetotal
          \advance\dimenone by#1         
          \setbox0=\vbox{#2\parskip=0pt                
                     \hyphenpenalty=10000
                     \rightskip=0pt plus 5em
                     \noindent#3 \vskip#4}    
        \advance\dimenone by\ht0
        \advance\dimenone by 3\baselineskip   
        \ifdim\dimenone>\pagegoal\vfill\eject\fi
          \else\eject\fi\fi}


\def\lsubheadline #1 #2{\bigbreak\vskip-\lastskip
      \checkleftspace{0.9cm}{\bf}{#1}{\bigskipamount}
         \vbox{\vskip0.7cm}\cen{\bf #1}\msk \cen{\bf #2}\bsk}

\def\sectionheadline #1{\bigbreak\vskip-\lastskip
      \checkleftspace{1.1cm}{\bf}{#1}{\bigskipamount}
         \vbox{\vskip1.1cm}\cen{\bfone #1}\bsk}

\def\lsectionheadline #1 #2{\bigbreak\vskip-\lastskip
      \checkleftspace{1.1cm}{\bf}{#1}{\bigskipamount}
         \vbox{\vskip1.1cm}\cen{\bfone #1}\msk \cen{\bfone #2}\bsk}

\def\lchapterheadline #1 #2{\bigbreak\vskip-\lastskip\indent\vskip3cm
                       \cen{\bftwo #1} \msk \cen{\bftwo #2} \gsk}
\def\llsectionheadline #1 #2 #3{\bigbreak\vskip-\lastskip\indent\vskip1.8cm
\cen{\bfone #1} \msk \cen{\bfone #2} \msk \cen{\bfone #3} \nobreak\bsk\nobreak}


\newtoks\literat
\def\[#1 #2\par{\literat={#2\unskip.}%
\hbox{\vtop{\hsize=.15\hsize\nin [#1]\hfill}
\vtop{\hsize=.82\hsize\nin\the\literat}}\par
\vskip.3\baselineskip}

\def\references{
\sectionheadline{\bf References}
\frenchspacing

\entries\par}

\mathchardef\emptyset="001F 
\def\address{Author: \tt$\backslash$def$\backslash$address$\{$??$\}$}

\def\abstract #1{{\narrower\baselineskip=10pt{\noindent
\eightbf Abstract.\quad \eightrm #1 }
\bigskip}}

\def\addresstwo{}

\def\dlastpage{\par\vbox{\vskip1cm\nin
\line{
\vtop{\hsize=.5\hsize{\parindent=0pt\baselineskip=10pt\nin\address}}
\quad 
\vtop{\hsize=.42\hsize\nin{\parindent=0pt
\baselineskip=10pt\addresstwo}}
\hfill} }}

\def\Box #1 { \msk\par\nin 
\centerline{
\vbox{\offinterlineskip
\hrule
\hbox{\vrule\strut\hskip1ex\hfil{\smc#1}\hfill\hskip1ex}
\hrule}\vrule}\msk }

\def\adots{\mathinner{\mkern1mu\raise1pt\vbox{\kern7pt\hbox{.}}
                        \mkern2mu\raise4pt\hbox{.}
                        \mkern2mu\raise7pt\hbox{.}\mkern1mu}}


\def\author{Wolfgang Bertram, Karl-Hermann Neeb}

\def\address
{Wolfgang Bertram 

Institut Elie Cartan 

Facult\'e des Sciences, Universit\'e Nancy I

B.P. 239

F - 54506 Vand\oe uvre-l\`es-Nancy Cedex

France

bertram@iecn.u-nancy.fr

} 

\def\addresstwo 
{Karl-Hermann Neeb

Technische Universit\"at Darmstadt 

Schlossgartenstrasse 7

D-64289 Darmstadt 

Deutschland

neeb@mathematik.tu-darmstadt.de}

\def\rightheadline{\tenbf\folio\hfil{\rm\title}\hfil}
\def\leftheadline{\tenbf\folio\hfil{\rm\title}\hfil}

\def\PE{\mathop{\rm PE}\nolimits}
\def\TKK{\mathop{\rm TKK}\nolimits}
\def\Gl{\mathop{\rm GL}\nolimits}

\def\Str{\mathop{\rm Str}\nolimits}
\def\pr{\mathop{\rm pr}\nolimits}
\def\mod{\mathop{\rm mod}\nolimits}
\def\ider{\mathop{\rm ider}\nolimits}
\def\Gr{\mathop{\rm Gr}\nolimits}

\def\Idem{\mathop{\rm Idem}\nolimits}

\def\ug{{\uline\g}} 
\def\uG{\uline G} 

\def\X{{\frak X}}

\nin
{\obeylines \parindent 0pt }
\vskip2cm
\centerline{\bfone Projective completions of Jordan pairs}
\centerline{\bfone Part I. The generalized projective geometry of a Lie algebra}

\gsk
\centerline{\bf\author}
\vskip1.5cm \rm

\nin {\bf Abstract.}
We prove that the projective completion $(X^+,X^-)$ of the Jordan pair
$(\g_1,\g_{-1})$ corresponding to a $3$-graded Lie algebra
$\g=\g_1 \oplus \g_0 \oplus \g_{-1}$ can be realized inside the space
$\cal F$ of inner $3$-filtrations of $\g$ in such a way that the
remoteness relation on $X^+ \times X^-$ corresponds to transversality
of flags. This realization is used to give geometric proofs of structure
results which will be used in Part II of this work in order to define
on $X^+$ and $X^-$ the structure of a smooth manifold
(in arbitrary dimension and over general base fields or -rings). 

\bigskip
\nin {\bf Contents.}

1. Three-graded and three-filtered Lie algebras

2. Tangent bundle, structure bundle and the canonical kernel 

3. The Jordan theoretic formulation

4. Involutions, symmetric spaces, and Jordan triple systems

5. Self-dual geometries and Jordan algebras

6. Functorial properties

7. Central extensions of three-graded Lie algebras

8. Grassmannian geometries and associative structures

\bigskip
\nin {\bf AMS-classification:} 17C30, 17C37

\bigskip
\nin {\bf Key words:} Jordan pair, generalized projective geometry,
3-graded Lie algebra, projective elementary group, symmetric space,
projective completion

\sectionheadline{Introduction}

A basic construction in linear algebra permits to imbed
an affine space $V$ into a projective space $X$
as the complement of a ``hyperplane
at infinity'' -- let us assume here for simplicity that
everything is defined over a commutative field $\K$, so
that $X$ may be seen as the projective space $\P(W)$ with
$W \cong V \oplus \K$. In the real or complex case, if the
dimension is finite or if $V$ is e.g. a Banach space, 
the projective space $X$ is a smooth
manifold with $V$ as a typical chart domain. An atlas of 
$X$ is obtained by taking all affine parts of
$X$ (all complements of hyperplanes of $X$); as is well-known,
change of charts is then given by  rational and hence
differentiable expressions. Similar constructions are known
for other manifolds $X$ such as Grassmannians, spaces of
Lagrangians or conformal quadrics.

In the present work we will construct such manifolds in a 
very general context, in arbitrary dimension and over 
general base fields or -rings instead of $\R$ or $\C$.
The present and first part contains the algebraic theory,
and Part II ([BN03]) contains the analytic theory.
For the case of base fields other than $\R$ or $\C$, we use
in Part II
suitable concepts of differential calculus and of smooth manifolds
 developed in [BGN03]  which, in the case of 
locally convex real or complex model spaces -- in particular,
for Banach and Fr\'echet spaces --, agree with the
usual concepts (but work more generally for manifolds modeled
on any Hausdorff topological vector space).
The present Part I is of independent interest since indeed
a good deal of the above mentioned constructions is purely
algebraic and admits a beautiful Lie- and Jordan theoretic
interpretation.
Geometrically, we  work in the context of
{\it generalized projective geometries} (introduced in [Be02]),
and algebraically, in the context of {\it $3$-graded Lie
algebras} which in turn correspond to {\it Jordan pairs}
(however, the paper is self-contained,  and we assume only
basic knowledge of Lie-algebras).
As in the ordinary projective case, it is a purely algebraic
problem to define the
chart domains, to give the precise description of the intersection
of chart domains and and to find explicit formulas for
the transition maps between different charts. Once this
is established, differential calculus  can be applied in order
to show in Part II that these structures are
differentiable under some suitable and natural assumptions.
In this way we not only obtain e.g. Grassmannian manifolds, Lagrangian
manifolds or conformal quadrics in arbitrary dimension over
$\K=\R,\C,\Q_p,\ldots$, but also a wealth of {\it symmetric spaces
(over $\K$)} which generalize
the {\it symmetric Banach manifolds} (see the
monograph [Up85]) but include many completely new examples that
had not been accessible before. 
The symmetric spaces thus constructed are precisely those which
are in the image of the {\it Jordan-Lie functor} (cf. [Be02],
[Be00]).

 Let us now describe the contents in some more detail.
Our basic objects are on the one hand {\it $3$-graded Lie algebras}, i.e.
Lie algebras of the form $\g=\g_1 \oplus \g_0 \oplus \g_{-1}$
satisfying the relations $[\g_\alpha,\g_\beta] \subset \g_{\alpha
+\beta}$, and on the other hand {\it $3$-filtered Lie algebras}, i.e.
Lie algebras $\g$ with a flag $\f: 0 = \f_2 \subset \f_1 \subset \f_0 \subset \g$ of subalgebras such
that $[\f_\alpha,\f_\beta] \subset \f_{\alpha+ \beta}$. For simplicity we shall also write 
these flags as pairs $\f= (\f_1,\f_0)$. 
If $\g$ is $3$-graded,
then $D(X)=i X$ ($X \in \g_i$) defines a derivation
of $\g$ such that $D^3=D$, 
called the {\it characteristic element}, and if $D$ is
inner, $D = \ad(E)$, $E$ will be called an {\it Euler operator}.
The {\it space of inner $3$-gradings of $\g$} is
$$
{\cal G} = \{ \ad(E) : \, E \in \g, \, \ad(E)^3  = \ad(E) \}.
$$
As usual in algebra, graded structures have underlying {\it filtered
structures}. However, for every $3$-grading, there are {\it two}
naturally
associated filtrations, $\f^+:=\f^+(D): \g_1 \subset \g_1 \oplus \g_0 \subset 
\g$ and $\f^-:=\f^-(D) = \f^+(-D):\g_{-1} \subset \g_{-1} \oplus \g_0 \subset \g$.
If 
$$
{\cal F} = \{ \f^+ (D) : \, D \in {\cal G} \}
$$
denotes the {\it space of inner $3$-filtrations of $\g$}, then we have
an injection
$${\cal G} \into {\cal F} \times {\cal F}, \quad 
D \mapsto (\f^+(D), \f^-(D)). $$

The spaces $\cal G$ and $\cal F$ carry many interesting geometric
structures; one may say that the pair $({\cal F} \times {\cal F},{\cal G})$
is a ``universal model of the generalized 
projective geometry associated to $\g$".
On ${\cal F} \times {\cal F}$ there is a natural
relation of being  transversal: two flags $\e = (\e_1, \e_0)$ and 
$\f = (\f_1, \f_0)$ are transversal if 
$$ \g = \e_1 \oplus \f_0= \f_1 \oplus \e_0. $$
Our key result on the
structure of $3$-graded Lie algebras (Theorem 1.6) affirms that 
${\cal G} \subset {\cal F} \times {\cal F}$ is exactly the set
of pairs of transversal inner $3$-filtrations of $\g$, and the
set $\f^\top$ of filtrations transversal to a given filtration
$\f$ carries canonically the structure of an affine space over
$\K$ with translation group $(\f_1,+)$. 
The {\it elementary projective group} $G=G(D)$ of the $3$-graded
Lie algebra $(\g,D)$ is the group of automorphisms of $\g$ 
generated by the
abelian groups $U^\pm = e^{\ad(\g_{\pm 1})}$; it acts on $\cal F$
and on $\cal G$. We realize the {\it projective completion} $(X^+,X^-)$
of the pair $(\g_1,\g_{-1})$ as the $G$-orbits in
$\cal F$ of the base points $\f^-$ and $\f^+$ 
such that $V^\pm := U^\pm.\f^\mp = (\f^\pm)^\top$ are 
``affine parts of $X^\pm$" (Theorem 1.12).
Summing up, the ``generalized projective geometry $(X^+,X^-)$''
is imbedded as a subgeometry in $({\cal F},{\cal F})$.

Using  this model, we have a natural definition
of the ``tangent bundle" $T{\cal F}$ of $\cal F$ and of a 
``structure bundle" $T'{\cal F}$ (taking the r\^ole of 
a cotangent bundle), and of sections of these bundles.
Thus we can define, in a purely algebraic context, 
{\it vector fields} on $\cal F$ as well as a certain operator between
$T'{\cal F}$ and $T{\cal F}$ called the {\it canonical kernel}
(Chapter 2). Over the affine parts $V^\pm$, the bundles and 
their sections can be trivialized, and it is seen that
our vector fields are actually quadratic polynomial and that the
canonical kernel coincides with the well-known {\it Bergman
operator} from Jordan theory (see below). 
Thus we get a very natural interpretation of 
the ``Koecher construction" which consists 
of realizing a $3$-graded Lie
algebra by quadratic polyomial vector fields (cf.\ also [Be00,
Ch.\ VII], where in the finite-dimensional real case
another natural interpretation of this construction
is given by using the {\it integrability of almost 
{{(para-)} complex} structures}).
This approach naturally leads to one of the main results to
be used in [BN03], namely the chart description of the action of 
$\Aut(\g)$ by ``fractional quadratic maps" (Theorem 2.8) .

In Chapter 3 we explain the link of the preceding results with
Jordan theory: the pair $(V^+,V^-)=(\g_1,\g_{-1})$ together
with the trilinear maps $T^\pm:V^\pm \times V^\mp \times V^\pm \to
V^\pm$ given by triple Lie brackets is a {\it (linear) Jordan pair},
and one can express in a straightforward way all relevant formulas from
the preceding chapter by these maps. Thus we obtain in a 
calculation-free way the Bergman-operator, the quasi-inverse
and many of their fundamental relations 
and thus get new and ``geometric" proofs of many Jordan theoretic
results.

In Chapter 4 we add a new structure feature, namely an
{\it involution} of the $3$-graded Lie algebra. It leads
to a bijection $p:X^+ \to X^-$ which is called a {\it polarity}
in case that there exist non-isotropic points $x$ (i.e., $p(x) \top x$).
Then the space of all non-isotropic points carries the
structure of {\it symmetric space over $\K$}.
We prove that the structure maps of this symmetric space 
are given by suitable Jordan-theoretic formulas (Theorem 4.4),
which will allow to conclude in Part II of this work that
these structure maps are differentiable and so we really deal with
symmetric spaces in the category of smooth manifolds.

Chapters 5 up to 8 contain further material that is not
strictly necessary for Part II of this work:
in Chapter 5 we discuss those geometries that correspond to
{\it unital Jordan algebras}: using our realization of $X^\pm$
as $G$-orbits in ${\cal F}$, they are characterized by the 
simple property that $V^+ \cap V^-$ is non-empty;
in particular, $X^+ = X^-$. An axiomatic characterization of the
``canonical identification of $X^+$ and $X^-$" has been given
in [Be03]; thanks to our model, things are considerably easier
here than in the axiomatic approach. 

In Chapter 6 some functorial aspects of our constructions are
investigated. It is shown that surjective homomorphisms of $3$-graded Lie 
algebras induce equivariant maps of the associated geometries and we also show that 
inclusions of inner $3$-graded subalgebras containing $\g_1 + \g_{-1}$ induce isomorphisms of 
the corresponding geometries. 

In Chapter 7 we discuss central extensions of inner $3$-graded Lie algebras. 
We show that for each central extension $q \: \hat \g \to \g$ of an inner $3$-graded 
Lie algebra $\g$ the extended Lie algebra $\hat\g$ carries a natural structure of an inner 
$3$-graded Lie algebra for which $q$ is a morphism of $3$-graded Lie algebras. We further 
construct a universal inner $3$-graded central extension of $\g$. 
We know from Chapter 6 that quotient maps induce maps on the level of geometries. 
For central extensions we show that these maps are isomorphisms. 

In the final Chapter 8, we look at an important class of geometries,
the {\it Grassmannian geometries}: let $R$ be an associative algebra
over the commutative ring $\K$, $V$ be a right $R$-module,
$\cal P$ the space of all $R$-linear projectors $V \to V$
and $\cal C$ be the space of all $R$-submodules of 
$V$ that admit a complement.
Then, by elementary linear algebra, 
the pair $({\cal C} \times {\cal C},{\cal P})$ 
has the main features of a generalized projective geometry
(Prop. 8.2, cf.\ also [Be01]), and in fact there is a homomorphism
into the geometry $({\cal F} \times {\cal F},{\cal G})$
with $\g=\gl_R(V)$ which induces isomorphisms on subgeometries
that are homogeneous under the elementary projective groups
(Theorem 8.4).  
Such geometries, called {\it Grassmannian geometries}, 
correspond to {\it special Jordan pairs}, i.e., to subpairs of
{\it associative pairs}. In particular, if $V=R$, then
the Grassmannian geometry can also be called the ``geometry
of right ideals of the associative algebra $R$"; it corresponds
to $R$, seen as a Jordan algebra over $\K$. 

\ssk
Finally, we would like to add some comments on related work
and on some open problems. 
The elementary projective group and the projective completion
of a general Jordan pair have been introduced by J. Faulkner 
([Fa83]), and results closely related to ours have been
obtained by O. Loos ([Lo95]). Their results are based on
the axiomatic theory of Jordan pairs ([Lo75]) and hence work
even for base rings in which $2$ or $3$ are not invertible. In contrast,
we work in a Lie theoretic context  and hence
assume throughout that $2$ and $3$ are invertible in $\K$.
However, it is possible to extend our approach also to the case
of a general base ring $\K$ -- see Remark 3.9.
Our results are more general in the sense that they apply to
general $3$-graded Lie algebras (not only to the Tits-Kantor-Koecher
algebra of a Jordan pair) and to the general automorphism group
$\Aut(\g)$ (and not only to the important special case 
given by transformations corresponding to quasi-inverses).
As a by-product, we get new proofs of many Jordan theoretic
results. It is in interesting open problem wether it is possible to
derive ``all" Jordan theoretic formulas in a similar geometric
way --  in particular, we would like to have in our model 
a ``geometric" proof of
the fundamental identity (PG2) of a generalized projective geometry
(cf.\ [Be02]) which is very closely related to
the {\it fundamental formula} of Jordan theory. 

Closely related results have also been obtained by Kaup ([Ka83])
and Upmeier ([Up85]) in the complex case in presence of a Jordan-Banach
structure. 
In fact, some arguments used to prove our Structure Theorem 1.6 have
been used by Kaup in the proof of his Riemann Mapping Theorem
(see the proof of  [Kau83, Prop. (2.14)] and the detailed version
of this in [Up85, Lemma 9.16]).  Our proofs are much simpler since 
we work directly with the $3$-graded Lie algebra, whereas
Kaup and Upmeier always use its homomorphic image 
realized by quadratic polynomial vector fields
(called {\it binary Lie algebras} in [Up85]).

The special case of Grassmannians, especially in the context of
Banach manifolds, has attracted much attention since it plays
an important r\^ole in differential geometetry  and is related
to several interesting differential equations -- see,
e.g., [D92], [DNS87]; our constructions are similar to, but
much more general than the ones described there.
For further references to constructions of manifolds in contexts
related to Jordan theory see Part II  ([BN03]);
cf. also [Io03] for an extensive bibliography.

\msk \nin {\bf Acknowledgement.}
The second named author
thanks the Institut Elie Cartan for its hospitality during
the time of prepration of a first draft of the present paper.

\msk \nin
{\bf Notation.} Throughout this paper,
$\K$ is a commutative ring with unit $1$
such that $2$ and $3$ are invertible in $\K$.
In Chapter 8, $R$ denotes a possibly non-commutative
ring which is a $\K$-algebra.

\sectionheadline{1. Three-graded and three-filtered  Lie algebras}

\nin {\bf 1.1. Three-graded Lie algebras.}
A {\it $3$-graded Lie algebra} (over $\K$) is a Lie algebra
over $\K$ of the form $\g=\g_{1} \oplus \g_0 \oplus \g_{-1}$
such that $[\g_k,\g_l] \subset \g_{k+l}$, i.e., 
$\g_{\pm 1}$ are abelian subalgebras which are $\g_0$-modules,
in the following often denoted by $V^{\pm}$ or $\g_\pm$,
and $[\g_1, \g_{-1}] \subset \g_0$.
The map $D:\g \to \g$ defined by $D X = iX$ for $X \in \g_i$
is a derivation of $\g$, called the {\it characteristic element
of the grading}. It satisfies the relation $(D-\id) D (D+\id)=0$,
i.e., $D^3 = D$; we say that it is a {\it tripotent derivation}.
Conversely, any tripotent derivation $D:\g \to \g$ is diagonizable
with possible eigenvalues $-1,0,1$ and corresponding decomposition
of $X \in \g$:
$$
X = X_1 + X_0 + X_{-1}, \quad \quad
X_0 = X - D^2 X, \quad
X_1 = {DX + D^2X \over 2}, \quad
X_{-1}={-DX +D^2 X \over 2}.
\eqno (1.1)
$$
Since $D$ is a derivation, this eigenspace decomposition is a
$3$-grading.
Therefore we may identify the {\it space of $3$-gradings of $\g$}
with the set 
$$
\tilde {\cal G} := \{ D \in \der(\g) : \, D^3 = D \}
$$
of tripotent derivations.  If $D=\ad(E)$ is an inner tripotent
derivation, then $E$ is called an {\it Euler operator}, and we
denote by
$$
{\cal G} :=\{ \ad(E) : \, E \in \g, \, \ad(E)^3 = \ad(E) \}
\eqno (1.2)
$$
the {\it space of inner $3$-gradings of $\g$}.
The {\it odd part} of the $3$-graded Lie algebra
$(\g,D)$ is $\g_{-1} \oplus \g_1$, and we say
that $(\g,D)$ is {\it minimal} if it is generated by its odd part, that is,
$\g_0$ is generated by the brackets $[\g_1,\g_{-1}]$.

The following degenerate cases may arise: $D^2 = \id$, then
$\g$ must be abelian, and we have merely a decomposition of a
$\K$-module into complementary subspaces;
$D^2 = D$, then $\g=\g_0 \oplus \g_1$ is nothing but a $\g_0$-module
$\g_1$, in particular, $D=0$ corresponds to the case $\g_1 = \{0\}$. 

\ssk
\nin
{\bf 1.2. The projective elementary group.}
Let $(\g,D)$ be a $3$-graded Lie algebra over $\K$.
For $x \in \g_{\pm 1}$, 
the operator  $e^{\ad x} = \1 + \ad x
+{1\over 2}(\ad x)^2$ is a well defined automorphism of $\g$.
(In order to check that $e^{\ad(x)}$ is an automorphism
we need that $\g$ has no $3$-torsion.) 
The group generated by these operators,
$$ 
G := G(D):= \PE(\g,D) := \la e^{\ad x} \: x \in \g_{\pm 1} 
\ra \subeq \Aut(\g),
$$
is called  the {\it projective elementary group of $(\g,D)$}  
(see Section 3.2 for the relation with the projective elementary
group defined in Jordan theoretic terms, as in [Fa83], [Lo95]). 
Sometimes it will be useful to have a matrix notation for elements
of $G$: if $g \in \Aut(\g)$, we let, with respect to the fixed
$3$-grading,
$$
g_{ij}:= \pr_i \circ g \circ \iota_j: \g_j \to \g_i, \quad \quad
i,j = 1,0,-1,
$$
where $\iota_j:\g_j \to \g$ are the inclusion maps and
$\pr_i \:=\pr_i(D): \g \to \g_i$ the projections, given by
$$
\pr_1= {D + D^2 \over 2}, \quad
\pr_0= \1-D^2, \quad \pr_{-1}= {D^2 - D \over 2},
\eqno (1.3)
$$
and write $g$ in ``matrix form"
$$
g =
 \pmatrix{g_{11} & g_{10} & g_{1,-1} \cr
g_{01} & g_{00} & g_{0,-1} \cr
g_{-1,1} & g_{-1,0} & g_{-1,-1} \cr}.
\eqno (1.4)
$$
The subgroups $U^\pm := U^\pm(D):= e^{\ad \g_{\pm}}$ of $G$ 
are abelian and generate~$G$. 
If the grading derivation is inner, $D=\ad(E)$, then
$$
\exp:  \g_{\pm 1}\to U^\pm, \quad 
X \mapsto e^{\ad(X)}
$$
is injective
since $v\in \g_{\pm}$ implies $e^{\ad v}.E = E\mp v$.
In the general case, we define the {\it automorphism group of $(\g,D)$} 
to be
$$
\Aut(\g,D) = \{ g \in \Aut(\g) : \, g \circ D = D \circ g \},
$$
and  
we further define subgroups $H:=H(D)$ and $P^\pm:=P^\pm(D)$ of $G$ via 
$$ 
H := G(D) \cap \Aut(\g,D)
\quad \hbox{ and } \quad 
P^\pm := H U^\pm = U^\pm H. 
\eqno (1.5)
$$
(If $D$ is inner, $D=\ad(E)$, then 
$ H= \{ h \in G \: h \circ \ad E \circ h^{-1}  = \ad E \} 
= \{ h \in G \: h.E - E \in \z(\g) \}$.)   
The groups $U^\pm$ are abelian, and since the group
$H$ commutes with $D$, it preserves the grading, hence 
normalizes $U^\pm$, so that $P^\pm$ are subgroups of $G$. 
Using notation from Equation (1.4),
the generators of $G$ are
represented by the following matrices (where $x \in \g_1$,
$y \in \g_{-1}$, $h \in H$):
$$
e^{\ad x}=
\pmatrix{\1 & \ad x & {1 \over 2} \ad(x)^2 \cr
0 & \1 & \ad x \cr 
0 & 0 & \1 \cr}, \ 
e^{\ad y}=\pmatrix{\1 & 0 & 0 \cr \ad y & \1 & 0 \cr 
{1 \over 2} \ad(y)^2 & \ad y& \1 \cr}, 
\  
h = \pmatrix{h_{11} & & \cr  & h_{00} & \cr  & & h_{-1,-1} \cr}.
$$
More information on the group $G(D)$ for inner $3$-gradings $D$
is given in Theorem 1.12.

Sometimes it will be useful to replace $G$ by a slightly bigger
group: if $D \in \tilde {\cal G}$ and 
$r \in \K^\times$, then, using the matrix notation (1.4),
$$
h^{(D,r)}  := \pmatrix{r &  & \cr & 1 & \cr & & r^{-1} \cr} =
r \pr_1 + \pr_0 + r^{-1} \pr_{-1} ,
\eqno (1.6)
$$
with the $\pr_i$ as in Equation (1.3),
defines an automorphism of $(\g,D)$ normalizing $U^\pm$
and commuting with all elements of the group $\Aut(\g,D)$. 
The group $G^{\rm ext}$ generated by $G$ and the group
$\{ h^{(D,r)} : \, r \in \K^\times \}$ will be
called the {\it extended projective elementary group}.

\ssk 
\nin
{\bf 1.3. Three-filtered Lie algebras.}
A {\it $3$-filtration} of a Lie algebra  $\g$
is a flag of subspaces 
$$
0=\f_{2 } \subset \f_{1} \subset \f_0 \subset \f_{-1}=\g
$$
such that
$$
[\f_k, \f_l] \subset \f_{k+l}.
\eqno (1.7)
$$
Supressing the trivial parts $\f_{2}$ and $\f_{-1}$ in the
notation, we
will denote such a flag by $\f=(\f_{1},\f_0)$ or 
$\f:(\f_1 \subset \f_0)$.
Let $\tilde{\cal F}$ be the set of such flags $\f$,
called  the {\it space of $3$-filtrations of
$\g$}. 
Conditions (1.7) are equivalent to the following requirements:

\ssk
\item{$\bullet$} $\f_0$ is a subalgebra, and $\f_{1}$ is an abelian
subalgebra of $\g$,
\item{$\bullet$} $\f_{1}$ is an ideal in $\f_0$, 
and $[\g,\f_{1}] \subset \f_0$.

\ssk
\nin
It follows that
the operators $\ad(X)$ with $X \in \f_{1}$ are 3-step nilpotent
and hence the automorphism $e^{\ad(X)}$ of $\g$ is well-defined.
We denote by
$$
U(\f):= e^{\ad(\f_{1})} =
\{ e^{\ad(X)} | \, X \in \f_{1} \} \subset \Aut(\g)
\eqno (1.8)
$$
the corresponding abelian group. 
{}From (1.7) it follows that $U(\f)$ preserves the filtration
$\f$. The filtration $\f$ is also stable under
the action of the subalgebra $\f_0$.

\msk \nin
{\bf 1.4. Relation between $3$-gradings and $3$-filtrations.}
To any $3$-grading $\g=\g_1 \oplus \g_0 \oplus \g_{-1}$ of $\g$
with characteristic derivation $D \in \tilde {\cal G}$
 we may associate {\it two} 
$3$-filtrations of $\g$, called the {\it associated 
plus- and minus-filtration},
given by the two flags 
$$
\f^+(D) := (\g_{1}, \g_0 \oplus \g_{1}), \quad \quad
\f^-(D) := (\g_{-1}, \g_0 \oplus \g_{-1}).
\eqno (1.9)
$$
Clearly, $\f^\pm(D) = \f^{\mp}(-D)$.
We will say that a $3$-filtration is {\it inner} if
it is of the form $\f = \f^+(\ad(E))=\f^-(\ad(-E))$ for some Euler operator
$E$, and the {\it space of inner $3$-filtrations}
will be denoted by
$$
{\cal F}:= \{ \f^+(D) : \, D \in {\cal G} \}.
\eqno (1.10)
$$
By these definitions, the maps
${\cal G} \to {\cal F}$, $D \mapsto \f^\pm(D)$ are surjective,
and the map
$$
{\cal G} \to {\cal F} \times {\cal F}, \quad
D \mapsto (\f^+(D),\f^-(D))
\eqno (1.11)
$$
is injective (since $\g_{\pm 1}$ are recovered by the filtration
and $\g_0 = (\g_0 \oplus \g_1) \cap (\g_{0} \oplus \g_{-1})$).

\ssk
\nin
{\bf 1.5. Transversality.} Two flags $\e = (\e_1, \e_0)$ and 
$\f = (\f_1, \f_0)$ as above are called {\it transversal} if 
$$ \g = \e_1 \oplus \f_0= \f_1 \oplus \e_0. $$
It is clear by construction
that the two filtrations $\f^+(D)$ and $\f^-(D)$ associated to
a $3$-grading $D$ of $\g$ are transversal.
We will prove that, conversely, any pair of transversal
inner $3$-filtrations arises in this way.
If $\f \in {\cal F}$, we will use the notation
$$
\f^\top := \{ \e \in {\cal F} : \, \e \top \f \}
\eqno (1.12)
$$
for the set of  inner $3$-filtrations that are transversal to $\f$,
and
$$
({\cal F} \times {\cal F})^\top = \{ (\e,\f) \in {\cal F} \times
{\cal F} : \, \e \top \f \}
\eqno (1.13)
$$
for the set of transversal pairs.

\Theorem 1.6. {\rm (Structure Theorem for the space of $3$-filtrations.)}
With the notation introduced above, the following holds for any Lie
algebra $\g$ over $\K$:
\item{(1)} The space of inner $3$-gradings can be canonically identified
with the space of transversal pairs of inner $3$-filtrations:
$$
{\cal G} = ({\cal F} \times {\cal F})^\top
$$
In other words,
two inner $3$-filtrations 
$\e$ and $\f$ are transversal if and only if
%
%
there exists an Euler operator $E$ such that
$\f =\f^+(\ad(E))$ and $\e = \f^-(\ad(E))$.
\item{(2)}
For any inner $3$-filtration $\f$, the space $\f^\top$ carries
a natural structure of an affine space over $\K$ with translation group $(\f_1, +)$. 
The group $\f_1$ acts simply transitively on $\f^\top$ by 
$x.\e := e^{\ad x}.\e$. 

\Proof. 
(1)
We have already remarked that ${\cal G} \subset ({\cal F} \times 
{\cal F})^\top$. In order to prover the other inclusion,
let us assume that $(\e,\f)$ is 
transversal. We have to show that $(\e,\f) \in {\cal G}$.

Since $\f$ is inner, there exists an Euler operator $E \in \g$ such that
$\f = \f^+(\ad(E))$. We also choose a $3$-grading 
$$ \g = \g_1 \oplus \g_0 \oplus \g_{-1} \quad \hbox{ with } \quad 
\e_1 = \g_1, \quad \e_0 = \g_1 \oplus \g_0. 
\eqno (1.14)
$$
Let $\pr_j:\g \to
\g_j$, $j=-1,0,1$, denote the corresponding projections. 
Our assumption that $\e$ and $\f$ are 
transversal means that
$$
\g= \e_1 \oplus \f_0 = \g_1 \oplus \f_0 
\quad\hbox{ and } 
\quad
\g= \f_1 \oplus \e_0 = \f_1 \oplus (\g_0 \oplus \g_1).
\eqno (1.15)
$$
{}From (1.14) together with
the first of these conditions we see that the restricted projection
$\pr_{-1}: \f_1  \to \g_{-1}$ is surjective.
Thus there exists a $Z \in \f_1$ such that
$\pr_{-1}(Z)=\pr_{-1}(E)$, i.e.
$Z - E \in \ker \pr_{-1} = \e_0$. 
Then $[Z,E]=-Z$ leads to 
$$
E':= e^{\ad(Z)} E = E + [Z,E] =
E - Z \in  \e_0. $$
On the other hand, since $Z \in \f_1$, the automorphism
$e^{\ad(Z)}$ stabilizes the flag $\f$, and hence
$$
\f^+(\ad(E'))= e^{\ad(Z)} \f^+(\ad(E))=
e^{\ad(Z)}.\f = \f.
$$
We may therefore, after replacing $E$ by $E'$, assume that $E \in \e_0$. 
This implies that $\e_1$ is invariant under $\ad E$, and since $\ad E$ acts by $-\id$ on the quotient 
space $\g/\f_0 \cong \e_1$, it follows that $[E,X] = -X$ for each $X \in \e_1$, hence 
$$ \e_1 =\{ X \in \g \: [E,X] = -X\}. $$
Further $\e_0 \oplus \f_1 = \g$, and $\f_1$ is the $1$-eigenspace of $\ad E$, so that 
the invariant subspace $\e_0$ must be the sum of the $0$- and $-1$-eigenspace of $\ad E$, and 
thus $\e = \f^-(\ad E)$. 


\ssk
(2) Assume $\e \top \f$. The group
$U(\f)= e^{\ad(\f_{1})}$
preserves both the flag $\f$ and the relation of being transversal.
 Therefore, for all $X \in \f_{1}$,
 $\f = e^{\ad(X)} \f$ is transversal to $e^{\ad(X)} \e$,
and hence we have an action of the abelian group $U(\f)\cong \f_1$
on $\f^\top$.

Let us prove that this action is transitive.
We assume that $\e' \top \f$ and $\e \top \f$.
By Part (1), there exists an Euler operator $E$ such that
$\f=\f^+(\ad(E))$ and 
$\e'=\f^-(\ad(E))$.
As in the proof of Part (1), we find $Y \in \f_1$ such that 
$$ E' := e^{\ad(Y)} E = E - Y \in \e_0. $$
Then the argument in (1) implies that $\ad E'$ is a characteristic element of the $3$-grading 
$$ \g_1 := \f_1, \quad \g_0 :=\f_0 \cap \e_0, \quad \g_{-1} := \e_1 $$
defined by the pair $(\f,\e)$. Therefore 
$$ \e = \f^-(\ad(E')) = e^{\ad Y}.\f^-(\ad(E)) = e^{\ad Y}.\e'. $$
Hence $\e'$ and $\e$ are conjugate under the group $U(\f)$.

Finally, we prove that the action is simply transitive.
Assume that $\e \in \f^\top$ and $X \in \f_1$ are such that
$e^{\ad(X)} \e = \e$.
Then $e^{\ad(X)}$ fixes the
transversal pair $(\e,\f)$ and hence commutes with
$\ad(E)$, where $E$ is an Euler operator such that $(\f,\e)=
(\f^+(\ad(E)),\f^-(\ad(E)))$. Applying this to
the element $E \in \g$, we get, since $[E,X]=-X$,
$$
0=e^{\ad(X)} \ad E(E) = \ad E \circ e^{\ad(X)}(E)=[E,E-X]=-X.
\qeddis


\Corollary 1.7.
Let $D_1=\ad(E_1), D_2=\ad(E_2) \in {\cal G}$ and
$\g_1 := \{ X \in \g | \, [E_1,X]=X \}$.
Then the following are equivalent:
\litem{(1)} $E_1$ and $E_2$ have the same associated $+$-filtration: 
$\f^+(\ad(E_1)) = \f^+(\ad(E_2))$. 
\litem{(2)} $D_1 - D_2 \in \ad(\g_1)$. 
\litem{(3)} $[D_1,D_2] = D_2 - D_1$. 
\litem{(4)} There is $v \in \g_1$ such that $D_2 = e^{\ad
v} D_1 e^{-\ad(v)}$. 

\Proof.
(4) implies (1) since $U(\f^+(D_1))$ preserves $\f^+(D_1)$.
Conversely, if (1) holds, then
$\f^-(D_2)$ is transversal to $\f^+(D_2)=\f^+(D_1)$,
and now (4) follows from  Part (2) of Theorem 1.6.

{}From (4) it follows that $\ad(E_2)=\ad(E_1 + [v,E_1] + {1 \over 2} [v,[v,E_1]])
=\ad(E_1-v)$, whence (2),
and from this it follows that
$[D_1,D_2]=\ad([E_1,E_2])=\ad([E_1,-v])=\ad(v)=\ad(E_1-E_2)$,
 whence (3), and finally
from (3) we get (4) by letting $v:=E_1-E_2$, which leads to 
$e^{\ad v}.E_1 = E_1 + [v,E_1] = E_1 + E_2 - E_1 = E_2$.
\qed

Next we state a ``matrix version" of Part (1) of Theorem 1.6, using
the matrix notation introduced in Equation (1.4).

\Corollary 1.8. With respect to a fixed inner $3$-grading
given by the Euler operator $E$, with corresponding pair
of $3$-filtrations 
$(\f^-,\f^+)= (\f^-(D), \f^+(D))= ((\g_{-1},\g_0 + \g_{-1}),(\g_1,\g_0 + \g_1))$, 
the following statements are equivalent:
\litem{(1)} $(g.\f^-,\f^+) \in {\cal G}$.
\litem{(2)} $\f^+$ and $g.\f^-$ are transversal. 
\litem{(3)}
$g_{-1,-1}$ and $(g^{-1})_{11}$ are invertible in $\End(\g_{-1})$, resp., 
in $\End(\g_1)$.

\Proof.
The equivalence of (1) and (2) is given by Theorem 1.6(1).
Now, (2) is equivalent to (4) and (5):

\litem{(4)} 
$g(\g_{-1})$ is a complement of $\g_{1} \oplus \g_0$ and
$g(\g_{-1} \oplus \g_0)$ is a complement of $\g_1$,
\litem{(5)}
$g(\g_{-1})$ is a complement of $\g_{1} \oplus \g_0$ and
$g^{-1}(\g_1)$ is a complement of $\g_{-1} \oplus \g_0$,

\ssk
\nin
and clearly (5) is equivalent to (3).
\qed

\Definition 1.9. For $x \in \g_1$ and $g \in \Aut(\g)$, 
we define
$$
d_g(x):=(e^{-\ad(x)} g^{-1})_{11}, \quad
c_g(x):=(g e^{\ad(x)})_{-1,-1}.
$$
Then
$$
d_g^+:=d_g:\g_1 \to \End(\g_1), \quad  \quad
c_g^+:=c_g:\g_1 \to \End(\g_{-1})
$$
are quadratic polynomial maps, called the {\it denominator}
and {\it co-denominator} of $g$ (w.r.t.\ the fixed inner grading
defined by $\ad(E)$). In a similar way, $d_g^-$ and $c_g^-$ are defined.
\qed

Writing $g$ and $e^{\ad(x)}$ in matrix form (1.4), the denominator for
$g^{-1}$ is given by
$$
d_{g^{-1}}(x)=
g_{11} - \ad(x) \circ g_{01} + {1 \over 2} \ad(x)^2 \circ g_{-1,1},
$$
and similarly for the co-denominator.
For the generators of $G$ we get  
the following (co-) denominators (where $v \in \g_1$, $w \in \g_{-1}$):
$$
\eqalign{
g=e^{\ad(v)}: & \quad d_g(x)= \id_{\g_1}, \quad c_g(x)=\id_{\g_{-1}} \cr
g=e^{\ad(w)}: & \quad
d_g(x)= \id_{\g_1} + \ad(x) \ad(w) + {\textstyle {1 \over 4}} \ad(x)^2 \ad(w)^2,  \cr
 & c_g(x) = \id_{\g_{-1}} + \ad(w) \ad(x) + 
{\textstyle {1 \over 4}}
\ad(w)^2 \ad(x)^2 \cr
g=h \in H: & \quad d_h(x)=(h_{11})^{-1} = (h^{-1})_{11}, \quad
c_h(x)=h_{-1,-1}. \cr}
\eqno (1.16)
$$
For $g=e^{\ad(w)}$ as in  the second equation, we introduce
the notation 
$$B_+(x,w):=d_g(x), \qquad B_-(w,x):=c_g(x). \eqno(1.17) $$
These linear maps define the {\it Bergman operator}, see Section 3.3.

\Corollary 1.10.
With respect to a fixed inner $3$-grading
given by the Euler operator $E$, we identify
$V^+:=\g_1$ with the set
$(\f^+)^\top = e^{\ad(V^+)} \f^-$.
Then for $x \in V^+$
the following statements are equivalent:
\litem{(1)} $(g.x,\f^+) \in {\cal G}$.
\litem{(2)} $\f^+$ and $g.x$ are transversal, i.e., $g.x \in V^+$. 
\litem{(3)}
$c_g(x)$ and $d_g(x)$ are invertible in $\End(\g_{-1})$, resp., in
$\End(\g_1)$.

\Proof.
This follows by applying Corollary 1.8 to the element
$g e^{\ad(x)} \in \Aut(\g)$.
\qed

In particular, for $g=e^{\ad(w)}$ with $w \in \g_{-1}$, 
it follows that $g.x \in V^+$ if and only if $B_+(x,w)$ and $B_-(w,x)$
are invertible.

\msk
\nin {\bf 1.11. The projective geometry of a $3$-graded Lie algebra.}
Recall from Section 1.2 the definition of the projective elementary 
group $G:=G(D)$.
Using Theorem 1.6(1), we may identify an
 inner grading $D=\ad(E)$ with the corresponding
pair $(\f,\e)=(\f^+(D),\f^-(D))$ of inner filtrations;
hence we may also write $G(\f,\e)$ for
the elementary group $G(D)$,
and similarly for $H(D)$ and $P^\pm(D)$.
If $\f,\e,\e'$ are inner $3$-filtrations such that
$\e \top \f$ and $\e' \top \f$,  
then Theorem 1.6(2) implies that $\e$ and $\e'$ are
conjugate under $G(\f,\e)$, and hence we have
$G(\f,\e)=G(\f,\e')$.
Therefore we may define the {\it projective elementary
group of the inner $3$-filtration $\f$} 
to be $G(\f):=G(\f,\e)$, where $\e \in {\cal F}$
is any filtration that is transversal to $\f$.
Note that
$$
U^+(\f,\e)=U^+(\f)
$$
is the abelian group defined by Eqn.\ (1.8) and hence is
independent of $\e$, whereas the groups
$U^-=U^-(\f,\e)$, $H=H(\f,\e)$ and $P^-=
P^-(\f,\e)$  depend on the choice of $\e$. (We will see below
that $P^+$ does not depend on $\e$.)
We define the following homogeneous spaces:
$$
X^\pm:=G/P^\mp, \quad \quad M:=G/H.
\eqno (1.18)
$$
 For reasons that will be
explained below, the data $(X^+,X^-,M)$ are called
the {\it (generalized) projective geometry associated
to the graded Lie algebra $(\g,D)$}.
The base point $(P^-,P^+)$ in $X^+ \times X^-$ will
often be denoted by $(o^+,o^-)$.

\Theorem 1.12.
 {\rm (Structure theorem for the projective
geometry of a $3$-graded Lie algebra)}
With the notation introduced above, the following holds:
\item{(1)}
The orbits of $D := \ad(E) \in {\cal G}$,
resp., of $\f^\pm \in {\cal F}$, under the action of $G$ are isomorphic
to $M$, resp., to $X^\pm$. In other words,
$$
H  = \{ g \in G(D): \, g.(\f^-,\f^+)=(\f^-,\f^+) \} \quad\hbox{ and } \quad  
P^\pm  = \{ g \in G(D) : \, g.\f^\pm  = \f^\pm \}. $$
Moreover, $P^+ \cap P^- = H$, $P^\pm \cap U^\mp = \{\1\}$ and
$$
P^\pm =\{ g \in G : g D g^{-1} - D \in \ad(\g_{\pm 1}) \}
=\{ g \in G \: g.E - E \in \z(\g) + \g_{\pm 1} \}.
$$
\item{(2)} If we identify $X^\pm$ with the corresponding orbits
in $\cal F$, then
$$
{\cal G} \cap (X^+ \times X^-) = M.
$$
\item{(3)}
For every element $\e \in X^-$, the set
$\e^\top$ is contained in $X^+$ and carries a well-defined
structure of an affine space over $\K$ with translation group $\e_1 = \g_1$.
In particular, $(o^-)^\top$ is canonically identified with
$V^+ = e^{\ad(\g_1)}. o^+$.
\item{(4)}
Consider the set $\Omega^+$ of elements of $G$ sending the base point 
$o^+ \in X^+$ to a point of the affine part $V^+ \subset X^+$,
$$
\Omega^+ := \{ g \in G : \, g.o^+ \in V^+ \}.
$$
Then the map
$$
\g_1 \times H \times \g_{-1} \to \Omega^+, \quad (v,h,w) \mapsto
e^{\ad(v)} h e^{\ad(w)}
$$
is a bijection, and moreover
$$
\Omega^+ = \{ g \in G : \, d_g(o^+) \in \Gl(\g_1), \,
c_g(o^+) \in \Gl(\g_{-1}) \}.
$$
\item{(5)}
The spaces $X^\pm \subset {\cal F}$ and $M \subset {\cal G}$
are stable under the action of
the automorphism group $\Aut(\g,D)$ and of the extended projective
elementary group $G^{\rm ext}$.

\Proof.
(1) An element $g \in G$ stabilizes $(\f^+,\f^-)$ if and only if it commutes with
$D=\ad(E)$ which means that it belongs to $H$. 

It is clear that $P^+$ stabilizes $\f^+$. Conversely,
assume that $g \in G$ satisfies $g.\f^+ = \f^+$.
Then $g.\f^+ = \f^+$ is transversal to $g.\f^-$, and hence by 
Theorem 1.6(2) there exists $v \in \g_1$ such that
$g.\f^- = e^{\ad(v)}\f^-$. Then $h:= e^{-\ad(v)} g$ preserves
$(\f^+,\f^-)$ and thus belongs to $H$. Therefore
$g = e^{\ad(v)} h$ belongs to $P^+$.
Hence $P^+$ is the stabilizer of $\f^+$. Similarly for $P^-$.

It follows that $P^+ \cap P^-$ is the stabilizer of $(\f^+,\f^-)$
which is $H$. Next, assume $g \in P^+ \cap U^-$.
Write $g = e^{\ad(v)}$ with $v \in \g_{-1}$. 
Since $v \mapsto e^{\ad(v)} \f^+$ is injective (Theorem 1.6(2)),
it follows from $g \f^+ =\f^+$ that $v =0$ and hence $g = \1$.

Finally, $g$ stabilizes $\f^+$ if and only if $D$ and $gDg^{-1}$ have
the same associated $+$-filtration, if and only if
$gDg^{-1} - D$ belongs
to $\ad(\g_1)$ (Corollary 1.7), whence the last claim of Part (1)
for $P^+$, and similarly for $P^-$.

(2) 
It is clear that the $G$-orbit $G.(\f^+,\f^-)$ belongs both to
$X^+ \times X^-$ and to $\cal G$. In order to prove the converse,
let $(\f,\e) \in (X^+ \times X^-) \cap {\cal G}$.
We may write $\f = g.\f^+$ for some $g \in G$.
Then $g^{-1}(\f,\e) = (\f^+,g^{-1} \e)$ again belongs to
$(X^+ \times X^-) \cap {\cal G}$.
According to Theorem 1.6, there exists $v \in \g_1$ such that
$g^{-1} \e = e^{\ad(v)} \f^-$.
It follows that $(\f,\e)=g e^{\ad(v)} (\f^+,\f^-)$
belongs to the $G$-orbit $G.(\f^+,\f^-)$.

(3) 
As in the proof of (2), we translate by an element
$g \in G$ such that $g \e = \f^-$, and then the claim
is precisely the one of Part (2) of Theorem 1.6.

(4)
Assume $g \in \Omega^+$ and let $v:=g.o^+ \in V^+$.
Then $e^{-\ad(v)} g.o^+ = o^+$, and according to Part (1), it follows that
then $p:=e^{-\ad(v)} g \in P^-$, whence the decomposition
$g = e^{\ad(v)}p=e^{\ad(v)}he^{\ad(w)}$.
Uniqueness follows from the fact that $P^+\cap P^- = H$.
Also, it is clear that any element $g \in U^+ P^-$
belongs to $\Omega^+$.

The second claim is a reformulation of Corollary 1.10.

(5) Assume $h \in \Aut(\g,D)$. From the relation
$h e^{\ad(x)} h^{-1}= e^{\ad(hx)}$ ($x \in \g_{\pm}$)
it follows that $h$ normalizes $G$.
Since $h$ stabilizes $\f^\pm$, it follows that, for
all $g \in G$, $hg.\f^-=hgh^{-1}.\f^- \in G.\f^- = X^+$.
It follows that $X^+,X^-$ and $M$ are stable under $h$.
Since $G^{\rm ext}$ is generated by $G$ and all $h^{(D,r)}$
(cf. Eqn. (1.6)), stability under $G^{\rm ext}$ also follows.
\qed

Since $P^+ \cap P^- = H$, 
$$
M \to X^+ \times X^-, \quad gH \mapsto (gP^-,gP^+)
\eqno (1.19)
$$
is a well-defined imbedding, and the following diagram commutes:
$$
\matrix{G/H & \into  & G/P^- \times G/P^+ \cr
\downarrow &  & \downarrow \cr
{\cal G} & \into & {\cal F} \times {\cal F}\cr}.
\eqno (1.20)
$$
Thus  we may say that the data $(X^+,X^-,M)$ forms
a subspace of $({\cal F},{\cal F},{\cal G})$ on which the
elementary projective group $G$ acts transitively.

\msk \nin
{\bf 1.13. The structure maps of the projective geometry.}
Assume $(\f_1,\f_2,\f_3)$ is a ``generic triple'' 
of inner $3$-filtrations;
by this we mean that it belongs to the space
$$
({\cal F} \times {\cal F} \times {\cal F})^\top :=
\{ (\f_1,\f_2,\f_3) \in
{\cal F} \times {\cal F} \times {\cal F} : \,
\f_1 \top \f_2, \, \f_3 \top \f_2 \}.
\eqno (1.21)
$$
Since $\f_2^\top$ carries a natural structure of an affine space over
$\K$ (Theorem 1.6(2)), we may take $\f_1$ as origin in 
$\f_2^\top$, i.e., we turn $\f_2^\top$ into a $\K$-module with
zero vector $\f_1$. Let $r \in \K$ and $r \f_3$ be the ordinary
multiple of $\f_3$ in this $\K$-module. Since it depends on $\f_1$ and
on $\f_2$, we write
$$
\mu_r(\f_1,\f_2,\f_3):= r_{\f_1,\f_2} (\f_3) := r \f_3 = (1-r) \f_1 + r \f_3, 
$$
where the latter expression only refers to the affine structure. 
The map 
$$
\mu_r:
({\cal F} \times {\cal F} \times {\cal F})^\top \to {\cal F}, \quad
(\f_1,\f_2,\f_3)  \mapsto \mu_r(\f_1,\f_2,\f_3)
\eqno (1.22)
$$
is called the {\it structure map of the projective geometry
$({\cal F},{\cal F},{\cal G})$}.
By restriction to the subgeometry $(X^+,X^-,M)$, we get in a
similar way two structure maps
$$
\mu_r^\pm:
 (X^\pm \times X^\mp \times X^\pm)^\top \to X^\pm $$
because for $\f_2 \in X^\mp$ we have 
$\f_2^\top \subeq X^\pm$ by Theorem 1.12(3). 
In [Be02, Th.\ 10.1] it is shown that these maps satisfy
two remarkable identities (PG1) and (PG2) which axiomatically
define the category of {\it generalized projective geometries}.
If $r \in \K^\times$, then we have
$$
\mu_r(\f_1,\f_2,\f_3) = h^{(D,r)} \cdot \f_3,
\eqno (1.23)
$$
where $h^{(D,r)}$ is the automorphism defined by Equation (1.6)
and $D$ corresponds to the 3-grading defined by the transversal
pair $(\f_1,\f_2)$.
The case $r=-1$ is of particular interest since it leads
 to associated {\it symmetric
spaces}, see Chapter 4.

\sectionheadline
{2. Tangent bundle, structure bundle and the canonical kernel}

\nin
{\bf 2.1. Tangent bundle and structure bundle.}
We continue to use the notation $\cal G$, resp.\ $\cal F$, for
the space of inner $3$-gradings (resp.\ $3$-filtrations) of
a Lie algebra $\g$.
For a $3$-filtration $\f=(\f_1,  \f_0)$,
we define $\K$-modules by
$$
T_\f {\cal F} := \g/\f_0, \quad
T_\f' {\cal F} := \f_1, 
\eqno (2.1)
$$
called the {\it tangent space of $\cal F$ at $\f$}, resp., the
{\it structural space of $\cal F$ at $\f$}.
If $\f=\f^-(\ad(E))$ is the minus-filtration w.r.t.\ an Euler operator
$E$, then $\f_0=\g_0 \oplus \g_{-1}$, and hence
$$
T_\f {\cal F} \cong \g_1, \quad \quad T_\f'{\cal F} = \g_{-1}.
$$
It is not misleading to think of $T_\f'{\cal F}$ as
a sort of ``cotangent space'' of $\cal F$ at $\f$. We let 
$$
T{\cal F}:= \bigcup_{\f \in {\cal F}} T_\f {\cal F},  \quad \quad
T'{\cal F}:=   \bigcup_{\f \in {\cal F}} T_\f' {\cal F}
\eqno (2.2)
$$
(disjoint union),
called the {\it tangent bundle} of $\cal F$, resp.,  the
{\it structure bundle} of $\cal F$.
The group $\Aut(\g)$ acts on $\cal G$ and on $\cal F$,
and  for any $g \in \Aut(G)$, the following maps are well-defined
and linear:
$$
\eqalign{
T_\f g:T_\f {\cal F} \to T_{g.\f} {\cal F}, & \quad Y \mod \f_0 \mapsto 
gY \mod g\f_0, \cr
T_\f'g:T_\f'{\cal F} \to T_{g.\f}'{\cal F}, & \quad Y \mapsto g Y, \cr}
\eqno (2.3)
$$
and if we define now $Tg:T{\cal F} \to T{\cal F}$, $T'g:T'{\cal F}
 \to T'{\cal F}$
in the obvious way, then  clearly the functorial
properties $T(g \circ h)=T(g) \circ T(h)$, and $T'(g \circ h)=T'(g) \circ
T'(h)$ hold.
Finally, if a base point $D \in {\cal G}$ is fixed and 
$X^\pm \subset {\cal F}$ are as in Cor.\ 1.10, then the tangent spaces
$T_\f X^\pm$, $T_\f' X^\pm$ and the corresponding bundles
$TX^\pm$, $T'X^\pm$ are defined.
The natural group acting on these spaces is the normalizer of
$G(D)$ in $\Aut(\g)$. 

\msk \nin
{\bf 2.2. Vector fields and the canonical kernel.}
If $Y \in \g$ and $\f \in {\cal F}$ is as above, we say that
$$
Y_\f := Y \mod \f_0 \in T_\f {\cal F}
\eqno (2.4)
$$
is {\it the value of $Y$ at $\f$}, and the assignment
$\tilde Y:{\cal F} \to T{\cal F}$,  $\f \mapsto Y_\f$ defines a
{\it vector field} on ${\cal F}$.
The space of vector fields on $\cal F$ is denoted by
$\X({\cal F})$; it is a $\K$-module in the obvious way such that
the surjection
$$
\g \to \X({\cal F}),\quad Y \mapsto \tilde Y
$$
becomes a $\K$-linear map which is equivariant w.r.t.\ the natural
actions of $\Aut(\g)$ on both spaces. In particular, the structural spaces
$T'_\f {\cal F}$ are subspaces of $\g$ and hence give rise to
vector fields. Composing with evaluation at another point, we are
lead to define, for
 $(\f,\e) \in {\cal F} \times {\cal F}$, a $\K$-linear map by
$$
K_{\f,\e}: T_\e' {\cal F}=\e_1  \to T_\f {\cal F}=\g/\f_0, \quad
Y \mapsto Y_\f = Y \mod \f_0.
\eqno (2.5)
$$
%
The collection of maps $(K_{\e,\f},K_{\f,\e})$, $(\f,\e) \in 
{\cal F} \times {\cal F}$,
is called the {\it canonical kernel}.
Note that $K_{\f,\e}$ is bijective if and only if $\e_1$ is a 
$\K$-module complement of $\f_0$ in $\g$. In particular,
if $\f=\f^-(\ad(E))$, $\e=\f^+(\ad(E))$, then
$K_{\f,\e}$ is identified with a linear map $\g_1 \to \g_1$
which is simply the identity.

\Theorem 2.3. For $\e,\f \in {\cal F}$ the 
following statements are equivalent:
\item{(1)} $(\e,\f) \in {\cal G}$,
\item{(2)} 
$K_{\f,\e}: T_\e' {\cal F} \to T_\f {\cal F}$ 
and $K_{\e,\f}:T_\f'{\cal F} \to T_\e {\cal F}$
are bijective.

\Proof.
The second condition clearly is equivalent to saying that
$\e$ and $\f$ are transversal, and therefore Theorem 2.3 is
a restatement of Part (1) of Theorem 1.6.
\qed

%
%

\msk \nin
{\bf 2.4. Trivialization over affine parts, and quadratic polynomial
vector fields.}
In the following we will often fix an Euler operator $E$, the associated
$3$-grading of $\g$ and the associated pair
$(\f^-,\f^+)=(\f^-(\ad(E)),\f^+(\ad(E)))$ of filtrations.
The pair $(\f^-,\f^+)$ then serves as a base point in $\cal G$
and in the homogeneous space $G.(\f^-, \f^+) \cong G/H \subset X^+ \times X^-$
(cf.\ Th. 1.12) and will also often be denoted by $(o^+,o^-)$.
The spaces $V^\pm:=\g_{\pm 1}$ are imbedded
into $X^\pm = G.\f^\mp \cong G/P^\mp$ via 
$X \mapsto e^{\ad(X)} \f^\mp$; this imbedding will
be considered as an inclusion, so that, for $x \in X^+$,
the condition $x \in V^+$ means that $(x,o^-) \in {\cal G}$.

%
%

The reader may think of $X^\pm$ as a kind of ``manifolds" modeled on
the $\K$-modules $V^\pm$: we will say that
$$
{\cal A}:= \{ (g(V^+),g): \, g \in G \}, \quad \quad
\phi_g:g(V^+) \to V^+, \quad g.x \mapsto x
\eqno (2.6) 
$$ 
is the {\it natural atlas of $X^+$}.
Having this in mind, a natural question is to describe the structures
introduced so far by a ``trivialized picture'' 
in the charts of the atlas $\cal A$.
Since the spaces $X^\pm$ are homogeneous under $G$, one can describe
$TX^\pm$ and $T'X^\pm$ as  {\it associated bundles}:
if $\pi:P^\pm \to \Gl(W)$ is a homomorphism of $P^\pm$ into the linear
group of a $\K$-module $W$, let
$$
G \times_{P^\pm} W = G \times W/ \sim
$$
with $(g,w) \sim (gp,\pi(p)^{-1}w)$ for $p \in P^\pm$.
If $\pi$ is the natural representation of $P^-$
on $W:=\g/(\g_0 \oplus \g_{-1}) \cong \g_1$ given by
$$
\pi(p):= p_{11}: \g_1 \to \g_1, \quad  X \mapsto \pr_1(pX)
\eqno (2.7)
$$
(this is the action of $P^-$ on $T_{\f^-} X^+$), then 
$$
G \times_{P^-} \g_1 \to TX^+, \quad
[g,X] \mapsto (T_{o^+} g) (X)
\eqno (2.8)
$$
is a $G$-equivariant bijection.
Similarly, if $\pi$ is the natural representation of $P^-$
on $W:=\g_{-1}$ given simply by $\pi(p)X=pX=p_{-1,-1}X$, then
$$
G \times_{P^-} \g_{-1} \to T'X^+, \quad
[g,X] \mapsto gX
\eqno (2.9)
$$
is a $G$-equivariant bijection.
For $TX^-$ and $T'X^-$ we have similar formulas.
If $f:G \to W$ is a function such that
$f(gp) = \pi(p)^{-1}.f(g)$ for all 
$g \in G$ and $p \in P^-$, then via
$$
s_f(g P^-) = [g,f(g)].
$$
we get a well defined section
of the natural projection $G \times_{P^-} W \to G/P^-$,
and every section arises in this way.
For instance, 
 for $Y \in \g$, the corresponding vector field
$\tilde Y$ on $X^+$ is given by the function
$$
\tilde Y_G:G \to \g_1, \quad
g \mapsto  g^{-1}Y \mod(\g_0 \oplus \g_{-1}) = \pr_1(g^{-1} Y),
\eqno (2.10)
$$
where  for the last equality we
 identified $\g/(\g_0 \oplus \g_{-1})$ and $\g_1$.
In fact, considering (2.8) as an identification, we have
$$
\tilde Y_{g.o^+} = Y \mod(g(\g_0 \oplus \g_{-1})) =
g(g^{-1}Y \mod(\g_0 \oplus \g_{-1})) = 
[g,g^{-1} Y \mod(\g_0 \oplus \g_{-1})] =
[g,\tilde Y_G(g)].
$$
We consider the special case $g=e^{\ad(v)}$ with $v \in \g_1$. 
We identify the restriction of $\tilde Y$ to $V^+ \subset X^+$
with the map
$$
\tilde Y^+ : V^+ \to V^+, \quad
v \mapsto  \pr_1(e^{-\ad v}.Y)
= \pr_1(Y-[v,Y]+{1 \over 2} [v,[v,Y]]).
\eqno (2.11)
$$
Note that the map $\tilde Y^+$ is a {\it quadratic} map
from $V^+$ to $V^+$.
In particular, it immediately follows from this formula
that for $Y \in \g_1$ this map is constant, for $Y \in \g_0$
it is linear and for $Y \in \g_{-1}$ it is homogeneous
quadratic:
$$ 
\tilde Y^+(v) = \cases{ 
Y & for $Y \in \g_1$, \cr
[Y,v] & for $Y \in \g_0$, \cr
{1\over 2} [v,[v,Y]] & for $Y \in \g_{-1}$. \cr} 
\eqno(2.12)
 $$
Similarly, $Y \in \g$ gives rise to a quadratic map
$\tilde Y^-:V^- \to V^-$.
Summing up, associating to $Y \in \g$ the quadratic  polynomial map
$\tilde Y^+ \times \tilde Y^- :V^+ \times V^- \to V^+ \times V^-$
 gives rise to
a {\it trivialization map}
$$
\g 
\to \Pol_2(V^+,V^+) \times \Pol_2(V^-,V^-)
$$
where $\Pol_2(W,W)$ is space of polynomial selfmappings of degree at most
two of a $\K$-module $w$.
Elements of $\g_0$ are mapped onto linear polynomials; in particular,
the Euler operator $E$ is mapped onto $(\id_{V^+}, -\id_{V^-})$.
The following result will not be used in the sequel, but is
recorded here for the sake of completeness.

\Proposition 2.5. 
The trivialization map is a homomorphism
of Lie algebras if we define the bracket of two quadratic
 polynomial maps
$p,q:W \to W$ on a $\K$-module $W$  by
$$
[p,q](x)=dp(x)q(x)-dq(x)p(x)
$$
where the (algebraic)
differentials $dp(x)$, $dq(x)$ of a (quadratic)
polynomial mapping
are defined in the usual way.

\Proof.
The commutator relations are directly checked by choosing
$p,q$ in the homogeneous parts $\g_1,\g_0,\g_{-1}$ of $\g$.
\qed

For the corresponding result on the group level, recall from
Definition 1.9 the nominator and co-denominator of an element
$g \in G$.

\Proposition 2.6.
If $g \in \Aut(\g)$ and $x \in V^+ \subset X^+$ are such that
$d_g(x)$ and $c_g(x)$ are invertible (equivalently, if $g.x \in V^+$),
then for all $Y \in \g$,
$$
\tilde{(g^{-1}Y)}^+(x) = d_g(x) \tilde Y^+ (g.x).
$$
In particular, for $Y=v \in \g_1$ we have
$$
\tilde{(g^{-1}v)}^+(x) = d_g(x) v.
$$
If $x, g_1.x$ and $g_1 g_2.x$ belong to $V^+$, then the
cocycle relation
$$
d_{g_1 g_2}(x)=d_{g_2}(x) \circ d_{g_1}(g_2.x)
$$
holds.

\Proof. The assumption that $g.x \in V^+$ means
that $g \circ e^{\ad(x)}$ belongs to the set $\Omega^+ \subset G$
defined in Theorem 1.12, Part (4).
Therefore, according to this theorem, 
 there exists a unique element $p(g,x) \in P^-$ such
that $g  e^{\ad(x)} =e^{\ad(g.x)} p(g,x)$ and hence
$p(g,x)=e^{-\ad(g.x)} g e^{\ad(x)}$. 
{}From this we get 
$$
\eqalign{
(p(g,x)^{-1})_{11} & =(e^{-\ad(x)} g^{-1} e^{\ad(g.x)})_{11} 
=\pr_1 \circ e^{-\ad(x)} g^{-1} e^{\ad(g.x)} \circ  \iota_1 \cr
& = \pr_1 \circ e^{-\ad(x)} g^{-1} \circ  \iota_1 
= (e^{-\ad(x)} g^{-1})_{11} = d_g(x). \cr}
$$
This will be used in the last line of the following calculation
(cf.\ also [Be00, VIII.B.2] for the general framework):
$$
\eqalign{
\tilde{(g^{-1}Y)}^+(x) & =
 \tilde Y_G (g e^{\ad(x)}) =  \tilde Y_G (e^{\ad(g.x)} p(g,x) ) \cr
&=  \pi(p(g,x))^{-1} \tilde Y_G(e^{\ad(g.x)}) 
=   (p(g,x)^{-1})_{11} \tilde  Y_G(e^{\ad(g.x)})  
=    d_g(x) \tilde Y^+(g.x). \cr}
$$
The second assertion follows since $\tilde v^+$ is a constant vector
field on $V^+$, see Equation (2.12).
The cocycle relation now follows:
$$
d_{g_1 g_2}(x)v=\widetilde{(g_2^{-1} g_1^{-1} v)}^+(x)=
d_{g_2}(x) \widetilde{ (g_1^{-1} v)}^+(g_2.x) =
d_{g_2}(x) \circ  d_{g_1}(g_2.x) v.
\qeddis

Proposition 2.6 implies in particular that 
the action of $g$ on the tangent bundle $TX^+$ is given in the canonical
 trivialization on $V^+$ by the expression
$T_x g \cdot v = d_g(x)^{-1} v$; in Part II of this work we
will show that, in presence of a differentiable structure,
this really corresponds to the differential $dg(x)$ of $g$
at $x$, applied to $v$. 
Similarly as in the proof of Prop. 2.6,
 it is seen that the action of $g$ on $T'X^+$ is, in
the trivialization $T'(V^+) \cong V^+ \times V^-$ over the affine
part $V^+ \subset X^+$, given by 
$$
T_x'g \cdot w = c_g(x) w,
$$
and that the co-denominators also satisfy a cocycle relation
$c_{g_1g_2}(x)=c_{g_1}(g_2.x) \circ c_{g_2}(x)$.

\msk
\nin {\bf 2.7. Nominators.}
We apply the preceding proposition in the case where $Y$ is
an Euler operator $E$  inducing the fixed $3$-grading of $\g$:
 for $g \in \Aut(\g)$ consider
the vector field $\tilde{g^{-1}E}$ on $X^+$ and define
the {\it nominator of $g$} to be the quadratic polynomial map
$$
n_g:V^+ \to V^+, \quad x \mapsto \tilde{g^{-1}.E}^+(x)=
\pr_1(e^{-\ad(x)} g^{-1} E)=(e^{-\ad(x)} g^{-1})_{10}.E. 
\eqno (2.13)
$$
Using the matrix notation (1.4), we can also write
$$
n_{g^{-1}} (x) = (g_{10} - \ad(x) \circ g_{00} + {1 \over 2}
\ad(x)^2 \circ g_{-1,0})(E).
$$
For the generators of $G$ we get the following nominators:
if $v \in \g_1$, $w \in \g_{-1}$, $h \in H$,
$$
n_g(x) = \cases{ 
x + v & for $g = e^{\ad(v)}$ \cr
x - {1 \over 2} \ad(x)^2 w & for $g = e^{\ad(w)}$ \cr 
x & for $g=h$ \cr} \eqno(2.14) 
$$
Note that the nominators will not depend on the Euler operator
$E$ such that $\ad(E)=D$ as long as $g$ acts trivially on the
center of $\g$; this is the case for all elements $g \in G$.
For general $g \in \Aut(\g)$ such that
$g.x \in V^+$, we can apply the preceding proposition and get,
using that $\tilde E^+(z)=z$ for all $z \in V^+$,
$$
n_g(x) = 
d_g(x) \tilde E^+(g.x) = d_g(x) (g.x).
$$
Since $d_g(x)$ is invertible, it follows that $g.x = d_g(x)^{-1} n_g(x)$.

\Theorem 2.8. Let $g \in \Aut(\g)$ and $x \in V^+$.
Then $g.x \in V^+$ if and only if
$d_g(x)$ and $c_g(x)$ are invertible, and then
the value $g.x \in V^+$ is given by
$$
g.x = d_g(x)^{-1} n_g(x).
$$
Using matrix notation $(1.4)$ and replacing $g$ by $g^{-1}$,
this can explicitly be written as an action of $\Aut(\g)$ on
$V^+$ by ``fractional quadratic maps'': if $g^{-1}.x \in V^+$,
then
$$
g^{-1}.x = 
(g_{11} - \ad(x) \circ g_{01} + {1 \over 2} \ad(x)^2 \circ g_{-1,1})^{-1}
\Big(g_{10} - \ad(x) \circ g_{00} + {1 \over 2}
\ad(x)^2 \circ g_{-1,0}\Big)(E).
$$

\Proof. For the first claim, see Corollary 1.10, and the second claim
is proved by the calculation preceding the statement of the theorem.
\qed

Using the formulas (1.16) for the denominators and (2.14) for the
nominators, we can now explicitly describe the fractional quadratic
action of the generators of $G$:
$$
\eqalign{
g=e^{\ad(v)}: & \quad g(x)=x+v \cr
g=e^{\ad(w)}: & \quad g(x)=\big(\id_{V^+} + 
\ad(x) \ad(w) + {\textstyle{1 \over 4}} \ad(x)^2 \ad(w)^2\big)^{-1}
(x-  {\textstyle{1 \over 2}} \ad(x)^2 w) \cr
g=h :  & \quad g(x) = h_{11} x. \cr}
$$

\msk
\nin {\bf 2.9. The automorphism group.}
The group $\Aut(\g,D)$ acts on $V^+ \times V^-$ by
$$
\eqalign{
\Aut(\g,D) & \to \Gl(V^+) \times \Gl(V^-),  \cr
h & \mapsto (h_{11},h_{-1,-1})=(d_{h^{-1}}(o^+),c_h(o^-))=
(d_h(o^+)^{-1}, c_h(o^-)). \cr}
$$
We denote by $\Aut_\g(V^+,V^-) \subset \Gl(V^+) \times \Gl(V^-)$
the image of this homomorphism (this is the automorphism group
of the associated Jordan pair, see Section 3.1), and by 
$\Str(V^+):=\pr_1 \circ \Aut_\g(V^+,V^-) \circ \iota_1$ its projection to
the first factor (sometimes called the {\it structure group
of $V^+$}). 

\Theorem 2.10.
If $x \in V^+$ and $g \in \Aut(g)$ satisfy $g.x \in V^+$, then
$d_g(x) \in \Str(V^+)$; more precisely,
$$
(d_g(x)^{-1},c_g(x)) \in \Aut_\g(V^+,V^-).
$$

\Proof.
If $g.x \in V^+$, then $g':=g e^{\ad(x)}$ belongs to the set
$\Omega^+ \subset G$ defined in Theorem 1.12.
According to Part (4) of this theorem, we decompose
$$
g'=e^{\ad(v)} h e^{\ad(w)}
 \eqno (2.15)
$$ 
with a unique $h=h(g,x) \in H$ depending
on $g$ and $x$.
{}From the definition of the (co-) denominators it follows then that
$$
d_g(x)=d_{g'}(0)=h_{11}^{-1}, \quad
c_g(x)=c_{g'}(0)=h_{-1,-1},
$$
and hence $(d_g(x)^{-1},c_g(x))=(h_{11},h_{-1,-1}) \in \Aut_\g(V^+,V^-)$.
\qed

As remarked after Proposition 2.6, the linear map $d_g(x)^{-1}$
can be intepreted as the tangent map of $g$ at $x$, and so
Theorem 2.10 means that $\Aut(\g)$ acts on $X^+$ by 
mappings that are {\it conformal with respect to the linear group
$\Str(V^+)$} (in the sense defined in [Be00, Section VIII.1.2]).
In some cases this already characterizes the group $\Aut(\g)$ 
as ``the conformal group of $X^+$''; this is the content of
the {\it Liouville theorem}, see [Be00, Ch. IX].

\sectionheadline
{3. The Jordan theoretic formulation}

\nin {\bf 3.1. Jordan pairs.}
If $(\g,D)$ is a 3-graded Lie algebra
and $V^\pm=\g_{\pm 1}$, the following $\K$-trilinear maps
are well-defined:
$$
\eqalign{
T^\pm : & V^\pm \times V^\mp \times V^\pm \to V^\pm, \cr
& (X,Y,Z) \mapsto T^\pm(X,Y,Z):= 
- [[X,Y],Z] = \ad(Z) \ad(X) Y 
= - \ad(X)\ad(Y) Z, \cr}
\eqno (3.1)
$$
and they satisfy the following identities,
where we use the notation $T^\pm(X,Y)Z:=T^\pm(X,Y,Z)$:
$$
\eqalign{
T^\pm(X,Y,Z) & = T^\pm(Z,Y,X), \cr
T^\pm(X,Y) T^\pm(U,V,W)&=
T^\pm(T^\pm(X,Y,U),V,W) \cr
& \quad \quad -T^\pm(U,T^\mp(Y,X,V),W) +
T^\pm(U,V,T^\pm(X,Y,W)),\cr}
\eqno (3.2)
$$
which mean that $((V^+,V^-),(T^+,T^-))$ is a
{\it linear Jordan pair over $\K$} (if $2$ and $3$ are 
invertible in $\K$, these two identities
imply  all other identities valid in Jordan pairs,
cf.\ [Lo75, Prop.\ 2.2(b)]). In the following we shall omit the
adjective {\it linear}, when dealing with Jordan pairs. 
Conversely, if $(V^\pm,T^\pm)$ is a
Jordan pair over $\K$, then
for $v \in V^\pm$ and $w \in V^\mp$ we define the operator 
$(v,w) \in \End(V^\pm)$ by 
$T^\pm(v,w).x := T^\pm(v,w,x)$ and 
let $\ider(V^+,V^-)\subeq \gl(V^+) \times \gl(V^-)$ 
be the Lie subalgebra generated by 
the operators $(-T^+(v,w), T^-(w,v))$, $v\in V^+$, $w \in V^-$.
The elements of this Lie algebra are called {\it inner derivations}. 
The algebra of {\it derivations of $(V^+,V^-)$} is defined by 
$$
\eqalign{
\der(V^+,V^-)&=\{ (A^+,A^-) \in \End_\K(V^+) \times \End_\K(V^+): \, 
(\forall u,v,w) \ 
\cr & \quad \quad 
A^\pm T^\pm(u,v,w)= T^\pm(A^\pm u,v,w)+T^\pm(u,A^\mp v,w)+
T^\pm(u,v,A^\pm w) \}, 
\cr}
\eqno (3.3)
$$
and it follows from (3.2) that it contains $\ider(V^+,V^-)$. 
Clearly, it contains also the element
$$
E := (\id_{V^+},-\id_{V^-}),
\eqno (3.4)
$$
called the {\it Euler operator of the Jordan pair $V^\pm$}.

If we are given a Jordan pair $(V^+, V^-)$, and 
$\g_0 \subeq \der(V^+,V^-)$ is a Lie subalgebra containing all
inner derivations, then there is a unique structure of
 a 3-graded Lie algebra on
$V^+ \oplus \g_0 \oplus V^-$
whose associated Jordan pair is $(V^-,V^+)$, and where the bracket
satisfies 
$$ 
[v,w] = (-T^+(v,w), T^-(w,v)), \quad v\in V^+, w \in V^-
\eqno (3.5)
$$
and the grading element is the Euler operator $E$ given by (3.4).
The subalgebra 
$$ 
\TKK(V^+,V^-) := V^+ \oplus (\ider(V^+,V^-) + \K E) \oplus V^- 
$$
is called the {\it Tits--Kantor--Koecher algebra of the Jordan pair 
$(V^+,V^-)$}. 
This choice for the $3$-graded Lie algebra associated to
$(V^+,V^-)$ has the advantage that $\z(\g)=0$.

The preceding construction may also be interpreted in the context
of {\it Lie triple systems} (cf.\ e.g. [Be00, Sect.\ III.3]):
it is essentially the {\it standard imbedding} of the (polarized)
Lie triple system $\q:=V^+ \oplus V^-$ into the corresponding
Lie algebra $\g=\q \oplus [\q,\q]$.
The standard imbedding yields a bijection between Lie triple
systems and Lie algebras with involution, generated by the 
$-1$-eigenspace of the involution.
See Chapter 6 concerning functorial properties of these constructions.

For any $\g_0$ as above, the representation of $\g_0$ on
$\g_{-1} \oplus \g_1$ will be faithful, so that 
$\z(\g) \cap \g_0 = \{0\}$.
It may happen for central extensions $\hat\g$ of $\g$ that
the corresponding subalgebra $\hat\g_0$ does not act faithfully on
$\hat\g_{-1} \oplus \hat\g_1 \cong \g_{-1} \oplus \g_1$ (see Chapter~7). 

\msk
\nin
{\bf 3.2. Projective elementary group and projective completion.}
For the rest of Chapter~3, we fix a Jordan pair $(V^+,V^-)$
and let $\g:=\TKK(V^+,V^-)$.
The projective elementary group $\PE(V^+,V^-):=G(\ad(E))$ 
is defined as in
Section 1.2. Using the notation, with $x,y \in V^\pm$, $v \in V^\mp$,
$$
\eqalign{
Q^\pm(x) v & := {1 \over 2} \ad(x)^2 v =
{1 \over 2} [x,[x,v]]  \cr
Q^\pm(x,y)  & := Q^\pm(x+y)-Q^\pm(x)-Q^\pm(y)= \ad(x)\ad(y):
V^\mp \to V^\pm, \cr}
\eqno (3.6)
$$
the  operators  $e^{\ad x} = \1 + \ad x
+{1\over 2}(\ad x)^2$ ($x \in V^\pm$)
are given in matrix notation by Equation (1.5), with
${1 \over 2} \ad(x)^2$ replaced by $Q^+(x)$ and
${1 \over 2} \ad(y)^2$ replaced by $Q^-(y)$. 
Our definition of $\PE(V^+,V^-)$ follows the one by O. Loos from [Lo95]. 
The projective linear group of a
Jordan pair has been introduced by Faulkner in [Fa83] in a slightly
different context (without Euler operator). 
The groups $P^\pm$ and the spaces $X^\pm = G/P^\mp$ are defined
as in Section 1.11; the embedding $V^+ \times V^- \to
X^+ \times X^-$ is called the {\it projective completion of the
Jordan pair $(V^+,V^-)$}.

\msk
\nin {\bf 3.3. The Bergman operator.}
Recall from Section 2.2 the canonical kernel: for $(x,y)\in
X^+ \times X^-$,
$$
K_{x,y}:T_y' X^- \to T_x X^+, \quad Y \mapsto Y_x.
$$
Of course, there is a similarly defined map $K_{y,x}$;
we will also use the notation $(K_{x,y}^+,K_{y,x}^-)$ for
$(K_{x,y},K_{y,x})$.
Using the description via associated bundles, the kernel is given
by
$$
 K_{g_1 P^-, g_2 P^+} \: T_{g_2P^+}' X^-  \to T_{g_1 P^-} X^+, \quad 
[g_2, v] \mapsto [g_1, \pr_1 (g_1^{-1}g_2.v)], 
\eqno(3.8)
 $$
and hence the trivialized picture is
$$
K^+_{x,y} = (e^{-\ad x} e^{\ad y})_{11} = d_{\exp -y}(x) \: V^+ \to V^+. 
\eqno (3.9)
$$
In matrix form,
$$
e^{-\ad x} e^{\ad y}  =
\pmatrix{\1 & -\ad(x) & Q^+(x) \cr
0 & \1 & -\ad(x) \cr 
0 & 0 & \1 \cr} \cdot 
\pmatrix{\1 & 0 & 0 \cr \ad(y) & \1 & 0 \cr Q^-(y) & \ad(y)& \1 \cr} ,
$$
so that we get for the coefficient with index $11$, using
that on $V^+$ we have for $x \in V^+$ and $y \in V^-$ the relation 
$\ad x \ad y = \ad [x,y] = - T^+(x,y)$:
$$
K^+_{x,-y} = B_+(x,y) = \id_{V^+}-T^+(x,y)+Q^+(x) Q^-(y). 
\eqno(3.10)
 $$
We likewise get 
$$ K^-_{x,-y} = B_-(y,x) =  \id_{V^-}-T^-(y,x)+Q^-(y) Q^+(x) $$
(cf.\ the definition in (1.17)). This expression 
is known as the {\it Bergman operator of the Jordan pair $(V^+,V^-)$}.
Theorem 2.3 now implies that 
{\it
the pair $(v,w)$ is transversal if and only if $(B_+(v,-w),B_-(-w,v))$ 
is invertible in $\End(V^+) \times \End(V^-)$.}
It is known in Jordan theory that $B_+(v,-w)$ is invertible if and only if so
is $B_-(-w,v)$ (the symmetry principle, cf.\ [Lo75, Prop. I.3.3]),
and hence $(v,w)$ is transversal if and only if $B_+(v,-w)$ is invertible.
So far we do not know a ``Lie theoretic'' proof of this fact.

\msk
\nin
{\bf 3.4. The quasi-inverse.}
Let $y \in \g_{-1}$. Then for $g = e^{\ad(y)}$ and $x \in V^+$,
Formulae (1.16) and  (2.14)  show that denominator, codenominator 
and nominator of $g$ are
given by
$$
d_g(x)= B_+(x,y), \quad \quad c_g(x) = B_-(y,x), 
\quad \quad n_g(x)= x - Q^+(x)y, 
\eqno (3.11)
$$
and hence, according to Theorem 2.8, 
$g(x) \in V^+$ if and only if $(B_+(x,y),B_-(y,x))$ is invertible, and then
$$
g.x =B_+(x,y)^{-1} (x - Q^+(x)y).
\eqno (3.12)
$$
Following [Lo77], we will use also the notation
$t_v (x) = x + v$ for translations on $V^+$ and 
$$
\tilde t_w(x):=e^{\ad(w)}.x= B_+(x,w)^{-1} (x - Q^+(x)w)
\eqno (3.13)
$$
for ``dual translations'' or ``quasi-inverses''.
 In Jordan theory  the notation $x^y := e^{\ad y}.x$
is also widely used (cf.\ [Lo75]), and one says that
$(x,y)$ is {\it quasi-invertible} if
$(B_+(x,y),B_-(y,x))$ is invertible, i.e., if
$(x,-y)$ is a transversal pair.
Our definitions of the Bergman operator via the canonical kernel
and of the quasi-inverse are natural in the sense that
they have natural transformation properties with respect
to elements $g$ of the group $\Aut(\g)$; taking for $g$ 
typical generators of $G$, we get Jordan theoretic results
such as the ``shifting principle'' (see [Be00, Section VIII.A]
for the precise form of the argument).

\msk \nin
{\bf 3.5. Automorphism and structure group.}
The group $\Aut_\g(V^+,V^-)$ defined in Section 2.9 coincides
for $\g = \TKK(V^+, V^-)$ with the automorphism group $\Aut(V^+,V^-)$ of $(V^+,V^-)$ in the 
Jordan theoretic sense. 
It follows from Theorem 2.10 that if $(x,-y)$ is transversal, then
$\beta(x,y):=(B_+(x,y),B_-(y,x)^{-1})$ belongs to $\Aut(V^+,V^-)$. 
The subgroup generated by these elements is called the {\it inner
automorphism group}. Projecting to the first factor, 
one gets the {\it structure group}, resp. the {\it inner
structure group of $V^+$}.

\msk \nin
{\bf 3.6.  Jordan fractional quadratic transformations.}
An {\it $\End(V^+)$-valued Jordan matrix coefficient
 (of type $(1,1)$, resp.\ of type $(1,0)$)}
is a map of the type
$$
q:V^\sigma \times V^\nu \to \End(V^+), \quad (x,y) \mapsto
(e^{\ad(x)}g e^{\ad(y)}h )_{11},
$$
where $\sigma,\nu \in \{ \pm  \}$ and $g,h$ belong to the
extended elementary projective group $G^{\rm ext}$ (cf.
Section 1.2), resp. 
$$
p:V^\sigma \times V^\nu \to V^+, \quad (x,y) \mapsto
(e^{\ad(x)}g e^{\ad(y)}h)_{10}E.
$$
These maps are quadratic polynomial in $x$ and in $y$.
Nominators and denominators of elements of $G$ are partial maps
of maps of the type of $p$ or $q$ by fixing one of the arguments to be zero. 
A {\it Jordan fractional quadratic map} is a
map of the form
$$
f:V^\sigma \times V^\nu \supset U \to V^+, \quad
(x,y) \mapsto q(x,y)^{-1} p(x,y),
$$
where $q,p$ are Jordan matrix coefficients of type (1,1),
resp.\ (1,0), and $U=\{ (x,y) \in V^\sigma \times V^\nu : \,
q(x,y) \in \Gl(V^+) \}$.
In the following, we also use the notation $\exp(x):=e^{\ad(x)}$
for $x \in V^\pm$.

\Theorem 3.7.
The actions
$$
V^+ \times X^+ \to X^+
\quad {\it and} \quad
V^- \times X^+ \to X^+
$$
are given, with respect to all charts from the atlas $\cal A$
(cf.\ {\rm Eqn.\ (2.6)}), by Jordan fractional quadratic maps.
In other words, for all $g,h \in G$, the maps
$$
(v,y) \mapsto (h \circ \exp(v) \circ g).y, \quad
(w,y) \mapsto (h \circ \exp(w) \circ g).y
$$
are Jordan fractional quadratic.

\Proof.  As to the first action, we write
$$
(h \circ \exp(v) \circ g).y=
(d_{h \circ \exp(v) \circ g}(y))^{-1}
n_{h \circ \exp(v) \circ g}(y) = q(v,y)^{-1} p(v,y)
$$
with 
$$
q(v,y) = d_{h \circ \exp(v) \circ g}(y)  =
(e^{-\ad(y)}g^{-1}e^{-\ad(v)} h^{-1})_{11}
$$
and
$$
p(v,y)= n_{h \circ \exp(v) \circ g}(y)
=(e^{-\ad(y)}g^{-1}e^{-\ad(v)}h^{-1})_{10}E,
$$
and hence the action is Jordan fractional quadratic.
For the action of $e^{\ad(V^-)}$, we use the same arguments.
\qed

We may say that $H_\infty:=X^+ \setminus V^+$ is the
``hyperplane at infinity''; then $H_\infty$ is stable
under the action of $V^+$. In case $(X^+,X^-)=(\K \P^n,
(\K \P^n)^*)$ is an ordinary projective geometry, the
action of the translation group on the hyperplane at
infinity is the trivial action. However, already in
the case of more general Grassmannian geometries this
is no longer true, as can be seen from the explicit formulas
for this case given in [Be01].


\Corollary 3.8.
With respect to the charts from the atlas $\cal A$, the
structure maps $\mu_r$ for $r\in \K^\times$
defined in {\rm Section 1.14} are given
by a composition of Jordan fractional quadratic maps
and diagonal maps $\delta(x)=(x,x)$.

\Proof. According to [Be02, Cor. 5.8], the multiplication maps
can be written as a composition of maps of the type described
in the preceding theorem, diagonal maps and one dilation $h^{(D,r)}$
(defined in Section 1.2). But
 this dilation comes from an element of $G^{\rm ext}$ and hence
the composition with such a dilation is  again Jordan fractional quadratic.
\qed

\msk
\nin {\bf 3.9. Case of a general base ring .}
Even if $\K$ is a general base ring (i.e. possibly 
with $2$ or $3$ not invertible), there still is a $3$-graded 
Lie algebra $\TKK(V^+,V^-)$ and a group $\PE(V^+,V^-)$ 
associated to a general (quadratic) Jordan pair, cf.\ [Lo95].  
The main difference is that 
in the matrix expression of $e^{\ad(x)}$ ($x \in \g_1$)
the term ${1 \over 2} \ad(x)^2$ has to be replaced by $Q^+(x)$.
Once one has checked that the abelian groups $U^\pm$ obtained
in this way are well-defined groups of automorphisms,
 one can essentially proceed as  
we did in Chapter 1, replacing the space ${\cal G}$
by the $\PE(V^+,V^-)$-orbit of $\ad(E)$ in $\der(\g)$ and
the space ${\cal F}$ by the space of inner filtrations 
belonging to gradings from $\cal G$.

\sectionheadline
{4. Involutions, symmetric spaces, and Jordan triple systems}

\msk \nin
{\bf 4.1. Symmetric spaces attached to a Lie algebra.} 
An {\it (abstract) reflection space} is a set $S$ together
with a map $\mu:S \times S \to S$ such that, if we let
$\sigma_x(y):=\mu(x,y)$,

\ssk
\item{(S1)} $\mu(x,x)=x$
\item{(S2)} $\sigma_x^2 =\id_S$
\item{(S3)} $\sigma_x$ is an automorphism of $\mu$, i.e.
$\sigma_x(\mu(y,z))=\mu(\sigma_x(y),\sigma_x(z))$.

\ssk
\nin (Differentiable reflection spaces, i.e. manifolds with a smooth
reflection space structure $\mu$,  have been introduced by
O. Loos in [Lo67]). In Part II ([BN03]) of this work we define
a {\it symmetric space (over $\K$)} to be a reflection space
$(S,\mu)$ such that  $S$ is a smooth
manifold over $\K$ (in the sense of [BGN03]) and $\mu$ is smooth
and satisfies

\msk
\item{(S4)}
the tangent map $T_x \sigma_x$ of $\sigma_x$ at $x$ is given
by $-\id_{T_x S}$.

\ssk
\nin  (See [BN03] for the basic theory of symmetric spaces
 and for a comparison with the approach by O. Loos [Lo69].)
To any Lie algebra $\g$ over $\K$ we may associate a reflection
space as follows. Let $S = \tilde{\cal G} = \{ D \in \der(\g) : D^3 = D\} $
be the space of $3$-gradings of $\g$ and
recall from Section 1.2 the definition of the extended
projective elementary group $G^{\rm ext}$ which is generated
by its normal subgroup $G$ and the subgroup
$\{ h^{(D,r)} | \, r \in \K^\times \}$.
Taking $r=-1$, we get the {\it reflection elements}
$$
\sigma^{(D)}:= h^{(D,-1)} = \1 -  2 D^2 \in \Aut(\g,D).
\eqno (4.1)
$$
We define the map $\mu$ by
$$
\mu: S \times S \to S, \quad
 \mu(D,D') := \sigma^{(D)} D' \sigma^{(D)} =
(\1-2D)D'(\1-2D).
\eqno (4.2)
 $$
Then (S1) follows from the fact that $D$ and $\sigma^{(D)}$ commute, 
(S2) holds because 
$\sigma^{(D)}$ is an involution, and (S3) follows from 
the fact that $\Aut(\g)$ clearly acts as a group of automorphisms of
$\mu$, and all reflection elements $\sigma^{(D)}$ belong
to $\Aut(\g)$.
It is clear that the subset ${\cal G} \subset \tilde{\cal G}$ 
is stable under $\mu$. 
Also, $M \subset {\cal G}$ is stable under $\mu$
because $M$ is stable under the action of $G^{\rm ext}$ 
(Theorem 1.12 (5)), and $G^{\rm ext}$
contains the reflection element $\sigma^{(D)}$ corresponding
to the base point and hence contains also all reflection elements 
corresponding to points of $M$.
Property (S4) is also satisfied in a purely algebraic sense:
since $\sigma^{(D)}$ acts by $-1$ on the complement
$\g_{\mp }$ of $\g_{\pm} \oplus \g_0$, it follows readily from the
definition of the tangent map in Section 2.1 that
$$
T_{\f^+(D)} \sigma^{(D)} = - \id_{T_{\f^+(D)}} {\cal F}, \quad
T_{\f^-(D)} \sigma^{(D)} = - \id_{T_{\f^-(D)}} {\cal F},
$$
and hence the tangent map $T_D(\sigma^{(D)})$ will be 
minus one if we define tangent map and tangent space at $D$
to be the direct product of the ones defined with respect
to $\f^+(D)$ and $\f^-(D)$.

The restriction of $\mu$ to ${\cal G} \times {\cal G}$ is related
to the ternary map $\mu_{-1}$ from Section 1.13 as follows:
assume $D_1$ corresponds to the transversal pair $(\f_1,\f_2)$ and
$D_2$ to the transversal pair $(\f_3,\f_4)$. Then
$$
\mu((\f_1,\f_2),(\f_3,\f_4)) 
= (\mu_{-1}(\f_1,\f_2,\f_3),\mu_{-1}(\f_1,\f_2,\f_4)),
= (h^{(D,-1)}.\f_3,h^{(D,-1)}.\f_4),
\eqno (4.3)
$$
which is the same as the product map on $M$ considered
in [Be02, Cor. 4.4].

\msk \nin
{\bf 4.2. Involutions and symmetric subspaces.}
An {\it involution} of a 3-graded Lie algebra is a Lie algebra
automorphism $\theta$ of order $2$ reversing the grading,
i.e., such that $\theta(\g_{\pm1})=\g_{\mp1}$
and $\theta(\g_0)=\g_0$. An involution $\theta$ 
induces by conjugation an automorphism of the elementary projective
group $G$, again denoted
by $\theta$, such that $\theta(P^-)=P^+$. Therefore it
induces a bijection
$$
p: X^+ \to X^-, \quad gP^- \mapsto \theta(g)P^+,
\eqno (4.4)
$$
compatible with the map ${\cal F} \to {\cal F}$, $\f \mapsto \theta(\f)$,
and
such  that $p(o^+)=o^-$.
We say that $\f \in {\cal F}$ is {\it non-isotropic (with respect to $\theta$)}
if $\theta(\f) \top \f$. In particular, the base point $o^+ = \f^-$
is non-isotropic; thus there exist non-isotropic points, and
$p$ is a {\it polarity} in the sense of [Be02].
Since $\theta$ is an automorphism normalizing $G$, the spaces
$\cal G$ and $M \subset {\cal G}$ are stable under $\theta$,
and the naturality of the product $\mu$ implies that
$\theta$ is an automorphism of $\mu$.
Therefore the $\theta$-fixed subspace $M^\theta$ is a symmetric
subspace of $M$, which as a set is in bijection with the set
of non-isotropic points of $X^+$, i.e.
$$
M^{(p)} := \{ \f \in X^+ : \, \f \, {\rm non-isotropic \,
w.r.t.} \, \theta \} \to M^{\theta}, \quad
\f \mapsto (\f,\theta)\f)
$$
is a bijection. By forward transport of structure, the
symmetric space structure of $M^\theta$ corresponds to
the structure on $M^{(p)}$ given by
$$
\mu(x,y)=\mu_{-1}(x,p(x),y)
\eqno (4.5)
$$
(this is the formula used in [Be02] to define the symmetric space
structure).
The symmetry w.r.t. the point $x$ is now
induced by the element $\sigma^{(D)}$, where $D \in {\cal G}$ corresponds to
the point
$(x,\theta(x)) \in {\cal G}$; as noticed above, the algebraically
defined tangent map
$T_x(\sigma^{(D)})$ equals minus the identity, and hence
(S4) is again satisfied in an algebraic sense.

\Theorem 4.3. 
For a fixed polarity $p:X^+ \to X^-$, we identify $X^+$
and $X^-$ via $p$. Then
the multiplication map $\mu$ on $M^{(p)}$ is a composition
of Jordan fractional quadratic maps and diagonal maps $\delta(x)=(x,x)$.

\Proof. 
By Corollary 3.8, the map $\mu_{-1}$ is of the form mentioned 
in the claim. According to Formula (4.5), $\mu$ is related to
$\mu_{-1}$ via
$$
\mu(x,y)=\mu_{-1}(x,x,y), \quad {\rm i.e.} \quad
\mu = \mu_{-1} \circ (\delta \times \id),
$$
which proves the claim.
\qed

In [BN03] it will be shown that Theorem 4.3 implies,
in very general situations, smoothness of $\mu$.

\msk
\nin {\bf 4.4.
Involutions and Jordan triple systems.} If $\theta$ is an involution
of the $3$-graded Lie algebra $\g$,
the trilinear map on $V^+$ defined by
$$
T(X,Y,Z):=-[[X,\theta(Y)],Z]
\eqno (4.6)
$$
is a {\it Jordan triple product}, i.e., it satisfies the identities
(3.1) with the superscripts $\pm$ omitted.
Conversely, given a {\it Jordan triple system over $\K$} (abbreviated
JTS)
(i.e., a $\K$-module with a $\K$-trilinear map satisfying the
above mentioned identities), we can define an involution on 
the Lie algebra $V^+ \oplus \der(V^+,V^-) \oplus V^-$ by
$$
\theta(v,(A,B),w) = (w,(B,A),v),
\eqno (4.7)
$$
and the associated JTS is the one we started with. In this way
we get a {\it bijection between Jordan triple systems over $\K$
and minimal $3$-graded Lie algebras with involution} (see Section 1.1).

\sectionheadline
{5. Self-dual geometries and Jordan algebras}

\msk
\nin {\bf 5.1. Self-dual geometries.}
We fix a $3$-graded Lie algebra $\g$ with grading induced by
the Euler operator $E$.
Recall our realization of $X^+$ and $X^-$ as $G$-orbits in
the space $\cal F$ of $3$-filtrations of $\g$.
Two cases can arise: 
either $X^+ \cap X^-$ is empty, or $X^+ = X^-$.
In the latter case we let $X:=X^+=X^-$, and again two cases
are possible: either

\ssk
\item{(a)}
$V^+ \cap V^-$ is empty, or
\item{(b)} $V^+ \cap V^-$ is not empty; 
then we say that the geometry given by $(\g,E)$ is {\it self-dual},
and we let
$V^\times := V^+ \cap V^-$.

\nin
An equivalent chararacterization of 
self-dual geometries is:
there are three points $\f_1,\f_2,\f_3 
\in X^+$ such that
$\f_1 \top \f_2$, $\f_2 \top \f_3$, $\f_3 \top \f_1$
(namely, take $(\f_1,\f_2)=(\f^-,\f^+)$ to be the base point and
$\f_3$ some element of $V^+ \cap V^-$).

\msk \nin
{\bf 5.2. The Jordan inverse.}
Assume that $(\g,E)$ is self-dual and fix some point
$\f \in V^+ \cap V^-$. We claim that there exists an involution $j$
of $\g$ (cf.\ Section 4.2) such that $j(V^+) \cap V^+ \not= \emptyset$.
In fact,  
let $W:=\f^\top$; then $W \subset X$ carries a natural 
structure of an affine space over $\K$ (Theorem 1.12(3)),
and by assumption  $o^+$ and $o^-$ belong to $W$. 
Let $e \in W$ be the midpoint
of $o^+$ and $o^-$ in the affine space $W$. 
Since $e \in W$, the pair $(e,f)$ is transversal and hence
corresponds to a $3$-grading $\g=\g_1' \oplus \g_0' \oplus \g_{-1}'$, i.e.
to an element $D' \in {\cal G}$. 
Let $j:= h^{(D',-1)} \in G^{\rm ext}$
be the automorphism that is minus one on $\g_1' \oplus \g_{-1}'$ 
and one on $\g_0'$.
Then $j$ fixes $(e,f)$ and
acts by the scalar minus one on the $\K$-module $W$ with zero
vector $e$. Since $e$ is the midpoint of $o^-$ and $o^+$, 
it follows that $j(o^-)=o^+$, and since obviously $j$ is of order 
two, it  is an involution. 
The condition $j(\f^-)=\f^+$ implies that $j(V^+) = j((\f^+)^\top)=(\f^-)^\top = V^-$. 
In particular, $V^+ \cap j(V^+)=
V^+ \cap V^- = V^\times$ is non-empty by assumption.
It contains the point $e=j(e)$. 

Now we apply Theorem 2.8 in order to derive an explicit formula
for $j$ in the chart $V^+$: 
for $v \in \g_1$, let ${\bf v}=jv \in \g_{-1}$; by Equation (2.12), 
$\bf v$ gives rise to the homogeneous quadratic vector field 
${\bf v}^+(x)=Q^+(x){\bf v}$ on $V^+$.
{}From Proposition 2.6 we now get
$$
d_j(x)v = (j^{-1}v)(x)= {\bf v}^+(x) =Q^+(x){\bf v}. 
$$
In a similar way we see that $c_j(x)=Q^-({\bf x})v$.
(In fact, since $j$ is an involution, 
$c_j(x)=j d_j(-x) j = j d_j(x) j$, and this  is invertible 
if and only if so is $d_j(x)$.)
Corollary 1.10 now shows that $j(x) \in V^+$ if and only if
$x$ belongs to the set
$$
V^\times = \{ x \in V^+ : \, Q^+(x) \, {\rm invertible} \, \}, 
$$
which is called  the {\it set of invertible elements in $V^+$}. 
The nominator of $j$ is
$n_j(x)=-x$ since $j$ reverses the grading (i.e.\ $jE + E \in \z(\g)$). Now
Theorem 2.8 shows that, for $x \in V^\times$,
$j(x) = - Q^+(x)^{-1} x$. The map 
$$
V^\times \to V^\times, \quad x \mapsto j(x):=- Q^+(x)^{-1} x
$$
is known as the {\it Jordan inverse}.

\msk \nin
{\bf 5.3. Jordan algebras.}
Notation being as above, note that by construction $e$ is
a fixed point of $j$, i.e., we have $e=j(e)=-Q(e)e$.
The tangent map $T_{e} j$ is $- \id_{T_{e} X^+}$.
Now 2.11 implies that also
$$
d_j(e) = T_{e} (-\id_W) = - \id_{V^+}.
$$
In conclusion, we have a Jordan triple system with an element
$e$ such that $Q(e)=-\id_V$.
It is is known that then 
$$
x \cdot y := -{1 \over 2} T(x,e,y)
$$
defines a Jordan algebra structure on $V$ with unit element $e$.
Conversely, every unital Jordan algebra arises in this way
(cf. [Lo75, I.1.10]).

%
%

\msk \nin
{\bf 5.4. The self dual geometry associated to a unital
Jordan algebra.}
Now assume that $(\g,E)$ is $3$-graded and
there exists $e^- \in \g_{-1}$ such that
$Q^-(e^-): V^+ \to V^-$ is  a bijection. 
Let $g:=e^{\ad(e^-)}$.
We claim that the flag $\f^+:\g_1 \subset \g_0 \oplus \g_{1}$ 
is transversal to the flag $g(\f^+)$:
first of all, for $v \in \g_1$,
$$
\pr_{-1}(g(v))=\pr_{-1}(v+[e^-,v]+Q^-(e^-) v)= Q^-(e^-)v,
$$
hence $\pr_{-1} \circ g \circ \iota_1$ is bijective and
thus $g(\g_1)$ is a complement of $\g_1 \oplus \g_0$.
Next, $g(\g_1 \oplus \g_0)$ is a complement of $\g_1$:
equivalently, $e^{-\ad(e^-)} \g_1$ is a complement of 
$\g_0 \oplus \g_1$, which is true by the same argument.
Hence, $g \f^+ \top \f^+$. With $(o^+,o^-)=(\f^-,\f^+)$, this
means that $g.o^- \in  V^+ \subset X^+$; but since $g \in G$,
this means that $X^- = X^+$.
Moreover, $o^- \in V^-$, and $e^{\ad(e^-)}$ acts as
a translation on $V^-$; therefore $g.o^- \in V^- \cap V^+$,
and it follows that the geometry is self-dual.
-- Summing up:

\Theorem 5.5. For a Lie algebra $\g$ with Euler operator $E$,
the following are equivalent:
\item{(1)} The geometry given by $(\g,E)$ is self-dual.
\item{(2)} There is an involution $j$ of $(\g,E)$ such that
$j(V^+) \cap V^+ \not= \emptyset$.
\item{(3)} The Jordan pair $(V^+,V^-)$ contains invertible elements.
\item{(4)} The Jordan pair $(V^+,V^-)$  comes from a unital Jordan
algebra $(V,E)$.

\Proof. (1) $\Rarrow$ (2) $\Rarrow$ (3) has been shown in 5.2, and 
(3) $\Rarrow$ (1) has been shown in 5.4. 
The equivalence of (3) and (4) is well-known (cf. [Lo75, I.1.10]; see 5.3). 
\qed

We do not know wether the condition 
 $X^+ = X^-$ alone already implies that 
$V^+ \cap V^- \not= \emptyset$ -- in the finite-dimensional case over
a field this certainly is true since then the ``hyperplane at infinity'' 
$X^+ \setminus V^+$ is an algebraic hypersurface, and hence
$V^+$ and $V^-$ must intersect if they are both included in $X^+$.
However, in infinite dimension the ``hyperplane at infinity'' may
become rather ``big'' and may very well contain some affine parts
 --  this problem is also discussed in  [Be03].

\sectionheadline
{6. Functorial properties}

\nin
{\bf 6.1. Functoriality problems.}
So far we have considered the following categories:
Jordan pairs $(V^+,V^-)$ over $\K$;
3-graded Lie algebras $(\g,D)$ over $\K$;
generalized projective geometries $(X^+,X^-)$ (these 
may be defined here simply as the geometries $(X^+,X^-)$
associated to a 3-graded Lie algebra);
associated reflection spaces $(M,\mu)$;
elementary projective groups $G=G(\g,D)$ associated to 
3-graded Lie algebras.
What are the functorial relations between these categories?
It is obvious that homomorphisms of $3$-graded Lie algebras
induce, by restriction to the pair $(\g_1,\g_{-1})$,
homomorphisms of Jordan pairs.
Other functorialily problems are less trivial:

\ssk
\item{(FP1)}
When does a homomorphism of Jordan pairs induce a homomorphism
of the associated Tits-Kantor-Koecher algebras?
\item{(FP2)}
When does a homomorphism of Jordan pairs induce a homomorphism
of the associated generalized projective geometries, resp.
of the associated reflection spaces?
\item{(FP3)}
When does a homomorphism of Tits-Kantor-Koecher algebras induce
a homomorphism of the associated elementary projective
groups?
\item{(FP4)}
When does a homomorphism of general 3-graded Lie algebras induce
a homomorphism of the associated elementary projective
groups?

\msk
\nin {\bf 6.2. Functoriality of the Tits-Kantor-Koecher algebra.}
In general, a homomorphism of Jordan pairs does not induce
a homomorphism of the associated Tits-Kantor-Koecher algebra.
In fact, as remarked in Section 3.1, the Tits-Kantor-Koecher
algebra $\TKK(V^+,V^-)$ may be seen as the standard imbedding of the
polarized Lie triple system $V^+ \oplus V^-$;
but the standard imbedding of a Lie triple system does in general
not depend functorially on the Lie triple system.
However, for {\it surjective} homomorphisms this is the case
(cf. [Lo95, Prop. 1.6]), and it is also true for
finite-dimensional semisimple Lie triple systems over
fields (cf. [Be00, Th. V.1.9]).

\msk \nin
{\bf 6.3. Functorialiy of the projective geometry and of the
reflection spaces.}
Any Jordan pair homomorphism $\phi^\pm:V^\pm \to (V')^\pm$ 
induces, in a functorial way, a well-defined map of geometries
$$
\eqalign{
\tilde \phi^\pm: &  X^\pm \to (X')^\pm,  \cr
e^{\ad(v_1)} e^{\ad(w_1)} \cdots e^{\ad(v_k)} e^{\ad(w_k)}.o^+ &  \mapsto 
e^{\ad(\phi^+(v_1))} e^{\ad(\phi^-(w_1))} \cdots
 e^{\ad(\phi^+(v_k))} e^{\ad(\phi^-(w_k))}.(o')^+, \cr}
$$
where $v_i \in V^+, w_i \in V^-$, $i=1,\ldots,k$, $k \in \N$
([Be02, Th.\ 10.1]); the main point here is that the geometry
$(X^+,X^-)$ can be described by generators (namely $(V^+,V^-)$) and relations
(with respect to the product maps $\mu_r$ from Section 1.13),
and Jordan pair homomorphisms are compatible with the relations.
(If the geometry is {\it stable} in the sense of [Lo95], then
these relations are  given by {\it projective equivalence},
 cf. [Lo77], [Lo95].) 
A homomorphism of geometries in the sense of [Be02] induces
a homomorphism of the corresponding reflection spaces
(because the reflection space structure is defined via the
maps $\mu_r$);
therefore Jordan pair homomorphisms always induce 
homomorphisms of associated reflection spaces.

In particular, an isomorphism of Jordan pairs induces a bijection
of geometries. Therefore, if two 3-graded Lie algebras have the
same Jordan pair $(\g_1,\g_{-1})$, then there is a canonical
bijection between the associated geometries.
(Cf. Th. 6.6 below for another, elementary proof.)
In particular, as long as we are only interested in the
associated geometry $(X^+,X^-)$ (e.g., in Part II of this
work) we may without loss
of generality assume that $\g$ is a Tits-Kantor-Koecher
algebra.

\msk
\nin
{\bf 6.4. Functoriality problem for the projective elementary group.}
Let $\phi:\g \to \g'$ be a morphism of $3$-graded Lie algebras.
One would like to define a homomorphism
$\tilde \phi: G \to G'$ of the associated elementary projective
groups by requiring that
$\tilde \phi(e^{\ad(v^\pm)}) = e^{\ad(\phi v^\pm)}$, but in general
this will not be well-defined. Therefore we introduce the group
$$
G(\phi):=
\{ g=(g_1,g_2) \in G \times G' : (\forall X \in \g)\ 
g_2 \phi(X)=\phi (g_1 X) \}.
$$
Then the projection $\pr_1:G(\phi) \to G$ onto the first factor is
surjective: in fact, the image of $\pr_1$ contains the generators of $G$
because all $g_1:=e^{\ad(x)}$, $x \in \g_\pm$, preserve the ideal
$\ker(\phi)$, and so with $g_2:=e^{\ad(\phi(x))}$ the pair
$(g_1,g_2)$ belongs to $G(\phi)$. 
Since $G$ is generated by $e^{\ad(\g_\pm)}$, it follows that
the projection $\pr_1$ is surjective.
The kernel of the projection $\pr_1$ is given by all elements of
the form $(\1,g_2)$ where
$g_2$ acts trivially on the subalgebra $\phi(\g) \subset \g'$.
Therefore, if $\phi$ is surjective, then $\pr_1$
is a bijection, and $\pr_2 \circ (\pr_1)^{-1}:G \to G'$
is the desired homorphism (see [Be00, Section I.3] for similar
considerations on the level of symmetric spaces).
Combining with 6.2, we see that surjective Jordan pair
homomorphisms induce (surjective) homomorphisms of associated
elementary projective groups (this result is also contained
in [Lo95, Prop. 1.6]).

The functoriality problem is now reduced to the case of injective  
homomorphisms.
In good cases, one may then hope to recognize $\pr_1: G(\phi) \to G$
as a sort of covering of $G$, and thus to view $\pr_2$ as a sort
of lift of the desired homomorphism to a covering group.

\msk\nin 
{\bf 6.5.  Problem (FP4) for isomorphisms of Jordan pairs.} 
 Let $\g$ be a $3$-graded Lie algebra $\g$ with grading element
$E$ and $\uline \g \subeq \g$ an inner $3$-graded subalgebra containing $\g_\pm$. We denote by $G$, resp. by $\uG$ the associated elementary
projective groups.
In the present section we will see that the injective homomorphism
$\uline \g \to \g$ (which induces an isomorphism of associated
Jordan pairs) induces a surjective homomorphism ``in the
opposite sense'': $G \to \uG$.
 In particular we shall give another and more elementary proof of the fact
that the associated homogeneous spaces are the same (cf. 6.3).
 As $\ug$ contains 
$\g_\pm$, it is invariant under the group $G$ generated by $e^{\ad \g_\pm}$. 
Moreover, $G$ acts trivially on the quotient space $\g/\ug$, because 
its generators  have this property, i.e., 
$g.x - x \in \ug$ for each $x \in \g$ and $g \in G$. 

\Theorem 6.6. There is a surjective restriction homomorphism  
$$ R \: G \to \uG, \quad g \mapsto g\res_{\ug} 
\quad \hbox{ with } \quad R^{-1}(\uline H) = H \quad \hbox{ and } \quad 
R^{-1}(\uline P^\pm) = P^\pm. $$
For the corresponding homogeneous spaces, we have 
$$ \uline G/\uline P^\pm \cong G/P^\pm 
\quad \hbox{ and } \quad 
 \uline G/\uline H \cong G/H $$
as homogeneous spaces of $G$. 

\Proof. First we observe that $R(U^\pm) = \uline U^\pm$ implies that $R$ is surjective.

Let $\ad_{\ug} \: \g \to \der(\ug)$ be given by 
$\ad_{\ug}(x) := \ad x\res_{\ug}$ and let $E$ be an Euler operator defining the grading of $\g$, resp., 
an Euler operator $\uline E \in \ug$ defining the grading on $\ug$. 
Then the ideal $\ker \ad_\ug$ of $\g$ is invariant under $\ad E$, hence adapted to the grading. 
For $x \in \g_\pm$ we have $\ad_\ug(x)(E') = [x,E'] = \mp x$, so that 
$$ \ker \ad_{\ug} \subeq \g_0, $$
and in particular $\ad_{\ug}$ is injective on $\g_+ + \g_-$. 
For $x = x_+ + x_0 + x_-$ with 
$x_\pm \in \g_\pm$ and $x_0 \in \g_0$ we have 
$$ [\ad_{\ug} E, \ad_\ug x] 
= \ad_{\ug} [E, x] = \ad_{\ug}(x_+ - x_-). $$
If this bracket vanishes, then $x_+ - x_- \in \ker \ad_{\ug} \subeq \g_0$
implies $x = x_0 \in \g_0$, i.e., we obtain the refined information 
$$ \ker \ad_\g \subeq \ad_\g^{-1}(\z_{\ad_\ug}(\ad_\ug E)) = \g_0. $$

Now let $g \in G$ with $R(g) \in \uline H$. 
For $x \in \g_0$ we then have 
$$ \ad_{\ug}(g.x) = R(g) \circ \ad_{\ug}(x) \circ R(g)^{-1}, $$
and all three factors on the right hand side commute with 
the grading derivation $\ad_{\ug} E$ of $\ug$. Hence 
$\ad_{\ug}(g.x)$ commutes with $\ad_{\ug} E$, and the argument from
above implies that $g.x \in \g_0$. On the other hand 
$R(g)$ preserves the grading of $\ug$, and hence in particular the subspaces 
$\g_\pm$. This means that $g$ preserves all eigenspaces of 
$\ad E$ on $\g$, and therefore that $g$ commutes with $\ad E$, 
so that $g \in H$. 
We conclude that 
$R^{-1}(\uline H) \subeq H,$
and the converse inclusion follows from the fact that the action of $H$ on $\ug$ preserves 
the grading $\ug = \g_+ \oplus (\ug \cap \g_0) \oplus \g_-$ of $\ug$. 

{}From $\uline P = \uline H \uline U^\pm$ and 
$R(U^\pm) = \uline U^\pm$, we obtain 
$$ R^{-1}(P^\pm) = R^{-1}(\uline H) U^\pm \subeq H U^\pm = P^\pm. $$
Since $R(P^\pm) = R(H) R(U^\pm) \subeq \uline H 
\uline U^\pm = \uline P^\pm$, the first 
assertion follows. 

For the homogeneous spaces, we now get 
$$ \uline G/\uline P^\pm 
\cong G/R^{-1}(\uline P^\pm) 
= G/P^\pm 
\quad \hbox{ and } \quad 
 \uline G/\uline H 
\cong G/R^{-1}(\uline H) 
= G/H. 
\qeddis

\sectionheadline{7. Central extensions of three-graded Lie algebras} 

\nin In this section $\K$ denotes a field with $2, 3 \in \K^\times$. 

\msk\nin {\bf 7.1.} Let $\g$ be a $3$-graded Lie algebra with grading element $E$. 
In this section we assume that $\g$ is generated by $E$ and $\g_\pm$,
i.e., that 
$$ \g_0 = \K E + [\g_+, \g_-].
\eqno(7.1) $$ 
We shall show that the homogeneous spaces associated of the 
elementary projective group of $\g$ do not change for central
extensions. Combining these results with those of the preceding
section, it follows that they only depend on the Jordan pair 
$(\g_+,\g_-)$. 

\Lemma 7.2. Let $q \: \hat\g \to \g$ be a central extension of $\g$, i.e., 
$q$ is surjective and $\ker q$ is a central subspace of $\hat\g$. 
We pick an element $\hat E \in \hat \g$ with 
$q(\hat E) = E$. Then $\ad\hat E$ is diagonalizable with the
eigenvalues $\{\pm 1,0\}$ and defines a $3$-grading 
$$ \hat\g = \hat\g_+ \oplus \hat\g_0 \oplus \g_- $$
such that $q$ is a morphism of $3$-graded Lie algebras. 

\Proof. First we observe that $q \circ \ad \hat E = \ad E \circ q$. From the relation 
$(\ad E)^3 = \ad E$ we derive that 
$$ 0 
= \big( (\ad E)^3  - \ad E\big) \circ q 
= q \circ \big( (\ad \hat E)^3  - \ad \hat E\big), $$
and hence that 
$$ \big((\ad \hat E)^3  - \ad \hat E\big)(\hat\g) \subeq \ker q \subeq
\z(\hat\g). $$
Applying $\ad \hat E$, we see that 
$$ (\ad \hat E)^4 = (\ad \hat E)^2, $$
i.e., 
$$ (\ad \hat E)^2 (\ad \hat E - \1) (\ad \hat E + \1) = 0. $$
Let 
$$ \hat \g = \hat\g_+ \oplus \hat\g_0 \oplus \hat\g_{-1} $$
be the generalized eigenspace decomposition for $\ad \hat E$. 
Then 
$$ \ad \hat E\res_{\hat\g_\pm} = \pm \id_{\hat\g_{\pm}} 
\quad \hbox{ and } \quad 
(\ad \hat E)^2.\hat\g_0 =\{0\}. $$
{}From $\ker q \subeq \z(\hat\g) \subeq \hat\g_0$, we derive that 
$q\res_{\hat\g_\pm}$ is injective and maps 
$\hat\g_\pm$ bijectively onto $\g_\pm$. Therefore 
$\g_0 = \K E + [\g_+, \g_-]$ leads to 
$$ \hat\g_0 = q^{-1}(\g_0) = \ker q + \K \hat E + [\hat\g_+, \hat
\g_-]. $$
As $[\hat\g_+, \hat\g_-] \subeq \ker \ad \hat E$, we conclude that 
$$ \hat \g_0 \subeq \ker \ad \hat E, $$
and hence that $\hat E$ is a grading element for the $3$-grading 
$\hat\g = \hat\g_+ \oplus \hat\g_0 \oplus \hat\g_-.$
\qed

\msk\nin {\bf 7.3.} If $\g$ is $3$-graded with grading element $E$ and 
$\z \subeq \g$ is a central subspace, then $\z \subeq \ker \ad E =
\g_0$, and the quotient map 
$q \: \g \to \g/\z$ is a central extension which is a morphism of
$3$-graded Lie algebras. 

This implies that for a central 
extension $q \: \hat\g \to \g$ for which 
$\hat\g$ is $3$-graded with grading element $\hat E$, the 
Lie algebra $\g$ is $3$-graded with grading element $E := q(\hat E)$, 
and Lemma 7.2 provides the converse information, that 
if $\g$ is $3$-graded with grading element $E$ 
and generated by $E$ and $\g_\pm$, then the Lie algebra 
$\hat\g$ has a natural $3$-grading defined by an element 
$\hat E$ with $q(\hat E) = E$ and $q$ is a morphism of $3$-graded Lie
algebra. Passing to the subalgebra generated by $\hat E$ and 
$\hat\g_\pm$, we even obtain a $3$-grading satisfying the same
condition as $\g$. In fact, 
$\h := \hat\g_+ + \hat\g_- + [\hat\g_+, \hat\g_-] + \K \hat E
\subeq \hat\g$ is a $3$-graded subalgebra
with $q(\h) = \g$, so that $\g \subeq \h + \ker q \subeq \h +
\z(\hat\g)$. In particular, $\h$ is an ideal of $\hat\g$.  

These consideration show that to understand central extensions of
$3$-graded Lie algebras, a natural context is given by those central
extensions $q \: \hat\g \to \g$ which are morphisms of $3$-graded Lie
algebras with grading element satisfying (7.1). 

\Lemma 7.4. Let $q \: \hat\g \to \g$ be a central extension of
$3$-graded Lie algebras with grading elements $\hat E$ and $E = q(\hat
E)$ satisfying {\rm(7.1)}. Then 
$$ q^{-1}(\z(\g)) = \z(\hat\g) $$
and therefore 
$\g/\z(\g) \cong \hat\g/\z(\hat\g).$

\Proof. Since $q$ is surjective, we have $\z(\hat\g) \subeq
q^{-1}(\z(\g))$. If, conversely, $q(x) \in \z(\g)$, then 
$[x, \hat\g] \subeq \ker q \subeq \z(\hat\g) \subeq \hat\g_0$. 
In particular, we obtain $[x,\hat E] \in \hat\g_0$ and therefore 
$x \in \hat\g_0$. This in turn implies $[x,\hat\g_\pm] \subeq
\hat\g_\pm$. As $q\res_{\hat\g_\pm}$ is injective, 
$[x,\hat\g_\pm] \subeq \ker q \cap \hat\g_\pm = \{0\}.$
Therefore $x$ commutes with $\hat\g_\pm$ and $\hat E$, hence is
central because $\hat\g$ is generated by $\hat E$ and $\hat\g_\pm$. 
\qed

\Corollary 7.5. If $\g$ satisfies {\rm(7.1)}, then $\z(\g/\z(\g)) =
\{0\}$.

\Proof. The adjoint representation 
$\ad \: \g \to \ad\g \cong \g/\z(\g)$ is a central extension
satisfying the assumptions of Lemma~6.3. Therefore 
$\ker \ad = \z(\g) = \ad^{-1}(\z(\ad\g))$
implies $\z(\ad\g) = \{0\}$. 
\qed

\Remark 7.6. (a) Let $\g$ a $3$-graded Lie algebra with grading element
$E$ and $\uline\g \trile \g$ the ideal $\ug$ generated by $E$ and $\g_\pm$ (see Section 6.5). 
We consider the Lie algebra homomorphism 
$$ \uline{\ad} \: \g \to \der(\uline\g), \quad 
x \mapsto \ad x\res_{\ug}. $$
In view of Corollary 7.5,  $(\uline{\ad})(\ug) \cong \ad \ug \cong \ug/\z(\g)$ is a center-free
$3$-graded Lie algebra satisfying (7.1). 

\par\nin (b) If $\g$ is a center-free $3$-graded Lie algebra satisfying
(7.1) and $(V^+,V^-) = (\g_+,\g_-)$ is the corresponding Jordan pair,
then the representation 
$$ \ad_{V^\pm} \: \g_0 \to \der(V^+,V^-), \quad 
x \mapsto (\ad x\res_{V^+}, \ad x\res_{V^-}) $$
is injective, and 
$$ \g \to V^+ \oplus \der(V^+,V^-) \oplus V^-, 
\quad x_+ + x_0 + x_- \mapsto 
(x_+, \ad_{V^\pm} x_0, x_-) $$
is an embedding of Lie algebras, where the right hand side carries the
bracket defined in Section~3.1. 

On the other hand, the subalgebra of $\g$ generated by
$\g_\pm$ is isomorphic to the corresponding subalgebra of 
$V^+ \oplus \der(V^+,V^-) \oplus V^-$, which is 
$\TKK(V^+,V^-)$. 
\qed

\Definition 7.7. Let $\g$ be a Lie algebra. We write 
$\la \g,\g\ra$ for the quotient of 
$\Lambda^2(\g)$ by the subspace generated by the elements of the form 
$$ [x,y]\wedge z + [y,z]\wedge x + [z,x]\wedge y,  $$
and write $\la x,y\ra$ for the image of $x \wedge y$ in $\la
\g,\g\ra$. 
Then $\la \g,\g\ra$ carries a natural Lie algebra structure satisfying
$$ [\la x,y\ra, \la x',y'\ra] = \la [x,y], [x',y']\ra, $$
and the map 
$$ b_\g \: \la \g,\g\ra \to \g, \quad \la x,y\ra \mapsto [x,y] $$
is a homomorphism of Lie algebras. 
\qed

\Theorem 7.8. Suppose that $\g$ is $3$-graded with grading element $E$
satisfying {\rm(7.1)}. If $\g$ is perfect, then we put $\tilde\g :=
\la \g,\g\ra$, and if $\g$ is not perfect, then we define 
$$ \tilde \g := \la \g,\g\ra \rtimes \K \tilde E, $$
where $\ad \tilde E$ satisfies 
$$ [\tilde E, \la x,y\ra] := \la [E,x],y\ra + \la x, [E,y]\ra = \la E,
[x,y]\ra. \leqno(7.2) $$
Then there is a unique Lie algebra homomorphism 
$$ q_\g \: \tilde\g \to \g
\quad \hbox{ with } \quad 
q_\g(\la x,y\ra) = [x,y] 
\quad \hbox{ and } \quad 
q_\g(\tilde E) = E. $$
This homomorphism is surjective with central kernel, hence a
central extension of $\g$. Moreover, it is weakly universal in the
sense that for any central extension 
$q \: \hat\g \to \g$ with a $3$-graded Lie algebra $\hat\g$ with
grading element $\hat E \in \hat\g$ there exists a unique 
Lie algebra homomorphism $\alpha \: \tilde\g \to \hat\g$ with 
$q \circ \alpha = q_\g$ and, if $\g$ is not perfect, with $\alpha(\tilde E) = \hat E$. 

\Proof. First we observe that $\g = [\g,\g] + \K E$. If $\g$ is perfect, then 
$$ b_\g \: \tilde\g := \la \g,\g\ra \to \g $$
is the universal central extension of $\g$. 
If $\g$ is not perfect, then $E \not\in [\g,\g] = \im(b_\g)$. 
Therefore $\g \cong [\g,\g] \rtimes \K E$. 

The Lie algebra $\der(\g)$ acts in a natural way by derivations on 
$\la \g,\g\ra$ via 
$$ d.\la x,y\ra = \la d.x, y \ra + \la x, d.y\ra. $$
We may therefore form the Lie algebra 
$\tilde \g := \la \g,\g\ra \rtimes \K \tilde E,$
where $\ad\tilde E$ satisfies (7.2). 

In both cases we obtain quotient homomorphisms 
$q_\g \: \tilde\g \to \g$ with $\ker q_\g = \ker b_\g \subeq
\z(\la\g,\g\ra)$. From (7.2) we derive that the action of $E$ on 
$\la \g,\g\ra$ annihilates $\ker b_\g$, so that 
$\ker q_\g$ is central in both cases. This means that 
$q_\g$ is a central extension, and Lemma 7.2 implies that 
$\tilde\g$ is $3$-graded with grading element $\tilde E$. 
Moreover, 
$$ \uline{\tilde\g} := \tilde\g_+ + \tilde\g_- + [\tilde\g_+, \tilde\g_-] + \K \tilde E $$ 
is an ideal of $\tilde\g$ with $\uline{\tilde\g} + \ker q_\g = \tilde\g$. 
In view of $\la E, E \ra = 0$, we have 
$$ [\tilde\g,\tilde\g] 
= \la [\g,\g], [\g,\g] \ra + \la E, [\g,\g]\ra 
= \la \g, [\g,\g] \ra + \la [\g,\g], E\ra 
= \la \g,\g\ra. $$
Therefore $[\tilde\g,\tilde\g] \subeq \uline{\tilde\g}$ implies 
$\tilde\g = \uline{\tilde\g}$ and hence that $\tilde\g$ satisfies (7.1). 

We claim that $q_\g$ is weakly universal as a central extension of $3$-graded
Lie algebras satisfying (7.1). So let 
$q \: \hat\g \to \g$ be a central extension. 
Then the bracket map $\hat\g \times \hat\g \to \hat\g$ factors through
an alternating bilinear map 
$$ b \: \g \times \g \to \hat\g \quad  \hbox{ with } \quad
b(q(x),q(y)) = [x,y], \quad x,y \in \hat\g. $$
Then the Jacobi identity in $\hat\g$ implies that 
$b$ satisfies the cocycle condition 
$$ b([x,y],z) +  b([y,z],x) +  b([z,x],y) = 0. $$
Hence there exists a unique linear map 
$$ \phi \: \la \g,\g \ra \to \hat\g 
\quad \hbox{ with } \quad \phi(\la x,y\ra) = b(x,y), $$ 
and it is easy to see that $\phi$ is a homomorphism of Lie algebras. 
Moreover, $\phi$ is a morphism of $3$-graded Lie algebras, because the
grading on $\hat\g$ is induced by the map 
$x \mapsto b(E,q(x))$. If $\g$ is not perfect, then 
$\hat\g$ is not perfect, and no grading element $\hat E \in \hat\g$ is
contained in $[\hat\g,\hat\g]$. We may therefore extend 
$\phi$ to a Lie algebra homomorphism 
$$ \phi \: \tilde\g \to \hat\g \quad \hbox{ with } \quad 
\phi(\tilde E) = \hat E. $$
This proves the weak universality of $\tilde\g$ as a $3$-graded Lie
algebra with grading element $\tilde E$. The map 
$\phi \: \tilde\g \to \hat\g$ is not uniquely determined by the
requirement that $q \circ \phi = q_\g$ because we may add 
any Lie algebra homomorphism $\psi \: \tilde\g \to \ker q$, which
corresponds to the ambiguity in the choice of the grading element 
$\hat E \in \hat\g$. Note that the commutator algebra of $\tilde\g$ is
a hyperplane, so that $\psi$ is determined by $\psi(\tilde E)$. 
\qed

\msk\nin {\bf 7.9. Central extensions have isomorphic geometries.}
Next we compare the groups 
$$ G \subeq \la e^{\ad \g_\pm} \ra \subeq \Aut(\g) 
\quad \hbox{ and } \quad 
\hat G \subeq \la e^{\ad \hat\g_\pm} \ra \subeq \Aut(\hat\g), $$
where $q \: \hat\g \to \g$ is a central extension of $3$-graded Lie
algebras satisfying (7.1). 
Since each element of $\hat G$ fixes the kernel $\z := \ker q$
pointwise, it induces an automorphism of $\g$, and we thus obtain a
group homomorphism 
$$ q_G \: \hat G \to G \quad \hbox{ with } \quad q_G(g) \circ q = q
\circ g, \ g \in \hat G, $$
because $e^{\ad \hat\g_\pm}$ is mapped onto $e^{\ad \g_\pm}$. 
The following theorem provides a short direct argument for the
isomorphy of the geometries associated to central extensions. Since
the corresponding Jordan pairs are the same, this could 
also be deduced from the general result mentioned in Section 6.2. 

\Theorem 7.10. We have 
$q_G^{-1}(H) = \hat H$ and 
$q_G^{-1}(P^\pm) = \hat P^\pm$. 
For the corresponding homogeneous spaces, we have 
$$ \hat G/\hat P^\pm \cong G/P^\pm 
\quad \hbox{ and } \quad 
 \hat G/\hat H \cong G/H $$
as homogeneous spaces of $\hat G$. 

\Proof. Since $q_G$ maps a generating subset of $\hat G$ onto a generating
subset of $G$, it is surjective. 

First we observe that for any $h \in \hat G$ we have 
$$ q_G(h).E - E = q(h.\hat E - \hat E). \leqno(7.3) $$

If $h \in \hat H$, then $h.\hat E - \hat E \in \z(\hat\g)$, and (7.3)
leads to 
$$q_G(h).E - E = q(h.\hat E - \hat E) \in q(\z(\hat\g)) = \z(\g),$$ 
and
hence $q_G(h) \in H$. Suppose, conversely, that 
$q_G(h) \in H$. Then (7.3) implies 
$$ h.\hat E - \hat E \in q^{-1}(\z(\g)) = \z(\hat\g), $$
so that $h \in \hat H$. 

Since $\hat P = \hat H \hat U^\pm$ and 
$q_G(\hat U^\pm) = U^\pm$, we have 
$$ q_G^{-1}(P^\pm) = q_G^{-1}(H) \hat U^\pm \subeq \hat H \hat U^\pm =
\hat P^\pm. $$
Further $q_G(\hat P^\pm) = q_G(\hat H) q_G(\hat U^\pm) \subeq H U^\pm
= P^\pm$, 
and we obtain $q_G^{-1}(P^\pm) = \hat P^\pm$. 

For the homogeneous spaces, we now get 
$$ G/P^\pm 
\cong \hat G/q_G^{-1}(\hat P^\pm) 
= \hat G/\hat P^\pm 
\quad \hbox{ and } \quad 
 G/H 
\cong \hat G/q_G^{-1}(\hat H) 
= \hat G/\hat H. 
\qeddis 

\Remark 7.11. We take a closer look at the kernel of $q_G$. 
Let $g \in \ker q_G \subeq \hat H$. Then 
$g$ preserves the grading of $\hat\g$. Since
$q\res_{\hat\g_\pm}$ is injective, we conclude that 
$g\res_{\hat\g_\pm} = \id_{\hat\g_\pm}$, and hence that 
$g - \id_{\hat\g}$ vanishes on the subalgebra generated by 
$\hat\g_\pm$. Moreover, $\im(g - \id_{\hat\g}) \subeq \ker q = \z$, so
that 
$$ g = \1 + D, $$
where $D \: \hat\g \to \z$ is a linear map. As $g$ is an automorphism,
it follows that $D \in \der(\hat\g)$, and hence that 
$[\hat\g,\hat\g] \subeq \ker D$. If $\hat\g$ is perfect, then $D$
vanishes, but if $\hat E \not\in [\hat\g,\hat\g]$, then 
$$ \Hom_{\rm Lie}(\hat\g,\z) = \Hom(\K \hat E, \z) \cong \z $$
describes the possibilities for $D$, which is determined by $D(\hat E)
\in \z$.

Since, for $h \in \hat H$ and $v \in \g_\pm$ we have 
$h e^{\ad v} h^{-1} = e^{\ad h.v}$, the 
condition $g\res_{\hat\g_\pm} = \id_{\hat\g_\pm}$ implies that 
$h$ commutes with the generating subset $e^{\ad \g_\pm}$, and hence
that 
$$ \ker q_G \subeq Z(\hat G). $$
This means that $q_G \: \hat G \to G$ is a central extension of
groups. 
\qed

\Example 7.12. We consider the case of a trivial Jordan pair 
$(V^+, V^-)$, i.e., all the maps $T^\pm$ vanish. Then the
corresponding $3$-graded Lie algebra is the semidirect sum 
$$ \g = (V^+ \oplus V^-) \rtimes \K E, $$
where 
$$ [E,(v_+, v_-)] = (v_+, - v_-) 
\quad \hbox{ and } \quad [V^+, V^-] = \{0\}. $$

Let $\beta \: V^+ \times V^- \to \z$ be any bilinear map. Then 
$$ \omega((v_+, v_-, \lambda E), (w_+, w_-, \mu E)) 
:= \beta(v_+, w_-) - \beta(v_-, w_+) $$
is a Lie algebra cocycle which defines a central extension 
$$ \hat \g = \g \oplus_\omega \z $$
with the bracket 
$$ [(x,z),(x',z')] := ([x,x'], \omega(x,x')), \quad x,x' \in \g, z,z'
\in \z. $$
The subalgebra of $\hat\g$ generated by $V^\pm$ is $2$-step nilpotent
and $\hat\g$ is solvable. In $\hat\g$ we have 
$$ [\ad V^+, \ad V^-] = \ad [V^+, V^-]\subeq \ad \z = \{0\}, $$
so that the groups $\hat G$ and $G$ are both abelian. 
Considering the orbit of the grading element, it is easy to see that 
$$ \hat G \cong V^+\times V^- \cong G. 
\qeddis 

\Remark 7.13. Let $\g$ be a $3$-graded Lie algebra with grading
element $E$. We have seen in Chaper~5 that the homogeneous spaces 
$G/H$ and $G/P^\pm$ are isomorphic to those associated to the 
subalgebra $\ug$ generated by $E$ and $\g_\pm$. 
Furthermore, the results in this section imply that the same holds for
the homogeneous spaces associated to the center-free Lie algebra 
$\ug/\z(\ug)$. The latter Lie algebra is isomorphic to the 
Tits--Kantor--Koecher Lie algebra 
$$\TKK(\g_+, \g_-) = \g_+ \oplus (\ider(\g_+,\g_-) + \K E) \oplus
\g_-$$ 
of the Jordan pair $(\g_+, \g_-)$. For that we only have to observe
that the triviality of the center implies that 
$\uline{\g}_0/\z(\g)$ embeds into $\der(\g_+,\g_-)$. 
We therefore obtain a natural
identification of the homogeneous space $G/H$ and $G/P^\pm$ with a
space of $3$-gradings of $\TKK(\g_+,\g_-)$, resp., 
a space of filtrations of this Lie algebra.
\qed


\sectionheadline
{8. Grassmannian geometries and associative structures}

\nin {\bf 8.1. Grassmannian geometries.}
Let $R$ be an associative algebra with unit $1$ over the
commutative unital ring $\K$ and let $V$ be a right
$R$-module. The {\it complemented Grassmannian (of $V$ over
$R$)} is the space
$$
{\cal C}:= \{ E \subset V : \, \exists F : V = E \oplus F \quad
(E,F: \, {\rm submodules \, of}\, V) \}
\eqno (8.1)
$$
of $R$-submodules of $V$ that have a complement.
For $V=R$ this is the {\it space of complemented right
ideals} of $R$ (cf.\ Section 8.6 below).
For $E,F \in {\cal C}$ we write $E \top F$ if 
$V = E \oplus F$; we let
$E^\top = \{ F \in {\cal C} : \, F \top E \}$ be
the set of complementary submodules of $E$ and
$$
({\cal C} \times {\cal C})^\top = \{ (E,F) \in {\cal C} \times
{\cal C}: \, V=E \oplus F \}.
\eqno (8.2)
$$
Let
$$
{\cal P} := \{ p \in \End_R(V) : \, p^2 = p \} = \Idem(\End_R(V))
\eqno (8.3)
$$
be the {\it space of projectors, resp., idempotents in $V$}. 
Taking $I:=2p - \id_V$ instead of $p$, we may also work with
the condition $I^2 = \id_V$ instead of $p^2 = p$ and view
$\cal P$ as the {\it space of polarizations of $V$}.
In this framework, the following analog of Theorem 1.6 is
an easy exercise in Linear Algebra:

\Proposition 8.2.
\item{(i)} The map
${\cal P} \to ({\cal C} \times {\cal C})^\top$, 
$p \mapsto (\im(p),\ker(p))=(\im(p),\im(\1-p))$ is a bijection.
\item{(ii)}
For all $E \in {\cal C}$, $E^\top$ carries canonically the
structure of an affine space over $\K$ (not over $R$
in general), modeled on the $\K$-module $\Hom_R(V/E,E)$.
\qed

\msk
\nin Moreover, $\cal P$ clearly is stable under
the binary map $\mu$ defined by $\mu(p,q)=(2p-\1)q(2p-\1)$
which defines on $\cal P$ the structure of a reflection
space (cf. 4.1). Using scalar extension by dual numbers
over $\K$, one may also define tangent bundles of
${\cal P}$ and ${\cal C}$, and then Property (S4)
will also hold for $\mu$; but we will not pursue this
construction here. 

\msk \nin {\bf 8.3. Flags and elementary group.}
We are going to describe the relation between this simple
linear algebra model and the model from Theorem 1.6.
Let $\g := \End_R(V)$ with the usual commutator as Lie
bracket. Note that the commutator is not $R$-bilinear
in general, but it is bilinear over the center of $R$;
hence $\g$ is a $\K$-Lie algebra.
An element $p \in {\cal P}$ defines a derivation 
$\ad(p)$ of $\g$ which is tripotent; with respect to the
decomposition $V=E \oplus F:=\im(p) \oplus \ker(p)$,
i.e., (in the obvious matrix notation)
$p=\bigl( {1 \atop 0} {0 \atop 0} \bigr)$, and
the grading of $\g$ is described by
$$
\eqalign{
\g_{-1} &  = \Big\{  \pmatrix{0 & 0 \cr \alpha & 0 \cr}
 : \, \alpha \in \Hom_R(E,F) \Big\}, \cr
\g_1 & = \Big\{  \pmatrix{0 & \beta \cr 0 & 0 \cr}:
 \, \beta \in \Hom_R(F,E) \Big\}, \cr
\g_0 & = \Big\{ \pmatrix{A & 0 \cr 0 & B \cr} :
\, A \in \End_R(E), B \in \End_R(F) \Big\}. \cr}
\eqno (8.4)
$$
Thus we have a well-defined map
from $\cal P$ to the space $\cal G$ of inner $3$-gradings of $\g$:
$$
\phi_{\cal P}: {\cal P} \to {\cal G}, \quad p \mapsto \ad(p).
\eqno (8.5)
$$
On the other hand, if $E \in {\cal C}$, then to the
``short flag" $0 \subset E \subset V$ we may associate a 
``long flag" $\f_E: 0 \subset \f_1 \subset \f_0 \subset \g$ by letting
$$
\f_1:=\{ X \in \g: \, X(V) \subset E, \, X(E)=0 \} \subset
\f_0:=\{ X \in \g: \, X(E) \subset E \} \subset \g;
\eqno (8.6)
$$
in matrix form: 
$$
\pmatrix{0 & * \cr 0 & 0 \cr} \subset
\pmatrix{* & * \cr 0 & * \cr} \subset
\pmatrix{* & * \cr * & * \cr}.
$$
It is clear that this is a $3$-filtration of $\g$
(even in an associative sense).
Thus we have a well-defined map
$$
\phi_{\cal C}: {\cal C} \to {\cal F}, \quad E \mapsto \f_E ,
\eqno (8.7)
$$
and it follows from the definitions that the diagram
$$
\matrix{
{\cal C} \times {\cal C} & \supset & {\cal P} \cr
\downarrow & & \downarrow \cr
{\cal F} \times {\cal F} & \supset & {\cal G} \cr}
\eqno (8.8)
$$
commutes. All maps in this diagram are obviously
equivariant with respect to the natural action of
the group $\Gl_R(V)$ on all spaces that are involved. 

If $E \in {\cal C}$ is fixed, then
the elements $X \in \f^1$ (with $\f^1$ as in (8.6)) are 
2-step nilpotent; thus $e^X = \1 + X$. Let
$$
U_E := e^{\f_1} = \1 + \f_1 =
\Big\{ \pmatrix{\1 & \beta \cr 0 & \1 \cr}
 : \, \beta \in \Hom_R(F,E) \Big\}, 
\eqno (8.9)
$$
where the latter matrix representation is with respect to a fixed
complement $F$ of $E$.
The group $U_E$ acts simply transitively on the set 
$E^\top$ of complements
of $E$. Therefore, if
for such a fixed decomposition $V=E \oplus F$, we let
$$
G(E,F):= \langle  U_E,U_F \rangle  \subset \Gl_R(V)
\eqno (8.10)
$$
be the group generated by $U_E$ and $U_F$, called the {\it elementary
group of $(V,E,F)$}, then $G(E,F)=G(E,F')$ for any two complements
$F,F'$ of $E$, and we may write also $G_E$ for $G(E,F)$.
We let
$$
P_E := \{ g \in G_E : \, g(E) =  E \}, 
\eqno (8.11)
$$
and, for a fixed complement $F$ of $E$,
$$
H(E,F) = \{ g \in G(E,F) : \, g(E)=E, \, g(F)=F \} = P_E \cap P_F.
\eqno (8.12)
$$

\Theorem 8.4. 
The equivariant maps $\phi_{\cal P}$ and $\phi_{\cal C}$ have
the following properties:
\item{(1)}
For all $E \in {\cal C}$, $\phi_{\cal C} (E^\top)=\phi_{\cal C}(E)^\top$.
\item{(2)}
For all $p \in {\cal P}$, the restriction of the map $\phi_{\cal P}$
to the orbit $\Gl_R(V).p$,
$$
{\cal P} \supset \Gl_R(V).p \to {\cal G}, \quad
g.p \mapsto \ad(g.p),
$$
is injective.
\item{(3)}
For all $E \in {\cal C}$, the restriction of the map $\phi_{\cal C}$
to the orbit $\Gl_R(V).E$,
$$
{\cal C} \supset \Gl_R(V).E \to {\cal F}, \quad
g.E \mapsto \f_{g.E},
$$
is injective.
\item{(4)}
Let $p \in {\cal P}$ with associated decomposition
$V=E \oplus F=\im(p) \oplus \ker(p)$.
The map $\phi_{\cal P}$ induces a bijection
$$
G(E,F)/H(E,F) \cong  G(E,F).p \to G(\ad(p)).\ad(p) \cong G(\ad(p))/H(\ad(p)),
$$
and the map $\phi_{\cal C}$ induces a bijection
$$
G(E,F)/P_F \cong  G(E,F) . E \to G(\ad(p)). \f_F = G(\ad(p))/P^-. 
$$

\Proof. 
(1) 
The action of $U_E$ on $\g$ is precisely the action
of $e^{\ad(\f_1)}$ on $\g$. 
Since $U_E$ acts simply transitively
on the set of complements of $E$, the claim follows from 
the corresponding fact about $\g$ ({\rm Theorem 1.6(2)}).

(2) 
Let $e \in {\cal P}$ and 
 $f := geg^{-1} \in {\cal P}$ with $g \in \Gl_R(V)$
such that $\ad(e)=\ad(f)$.
 Then $z:=f - e \in Z(A)$ where $A$ is the associative
$\K$-algebra $\End_R(V)$. In particular,
 $ef = fe$ and therefore 
$$ (e-f)ef = e^2 f - ef^2 = ef-ef=0.$$ 
We have 
$$ (f-e)^2 = f^2 - 2 ef + e^2 = f + e - 2ef $$
and 
$$ (f-e)^3 = (f-e)(f + e - 2ef) = f^2 - e^2 - 2(f-e)ef = f-e, $$
i.e., $z^3=z$.
Write $z=z_1-z_2$ with
$$ z_1 = {1\over 2}z(z+\1) \quad\hbox{ and } \quad 
z_2 = {1\over 2}z(z-\1). $$
Then $z_1$ and $z_2$ are again central, and $z_1^2=z_1$ and $z_2^2=z_2$.
This implies 
$$ z_1 
= {1\over 2}(f-e)(f-e+\1) 
= {1\over 2}(f+e - 2ef + f-e) 
= {1\over 2}(2f - 2ef) 
= f - ef = zf $$
and 
$$ z_2
= {1\over 2}(f-e)(f-e-\1) 
= {1\over 2}(f+e - 2ef - f+e) 
= {1\over 2}(2e - 2ef) 
= e - ef = -ze. $$
We further obtain 
$$  z_2 = - z_2 z = z_2(e - geg^{-1})  
= z_2 e - gz_2 eg^{-1} 
= -ze^2 + gze^2 g^{-1} 
= -ze + gzeg^{-1} 
= z_2 - g z_2 g^{-1} =0 $$
because $z_2$ is central, and likewise 
$$ z_1 = z_1 z = z_1(f - g^{-1}fg)  
= z_1 f - g^{-1} z_1 f g 
= z f^2 - g^{-1} z f^2 g 
= z f - g^{-1} z f g 
= z_1 - g^{-1} z_1 g = 0. $$
Eventually we obtain $z = z_1 - z_2 = 0$ and hence
$e=f$, as had to be shown. 

(3)
This follows by combining (2) and (1)
(observing that the fibers of the map ${\cal P} \to {\cal C}$,
$p \mapsto \im(p)$ are of the form $F^\perp$,
$F \in {\cal C}$, and similarly for ${\cal G} \to {\cal F}$).

(4)
This follows from (2) and (3), observing that
the action of $H(E,F)$ on $\g$ coincides with the action of $H(\ad(p)$. For 
$P_E$ we argue similarly. 
\qed


\msk
\nin {\bf 8.5. Special Jordan pairs.}
If $p \in {\cal P}$ and $\g=\gl_R(V)$, then the associated
Jordan pair is
$$
(\Hom_R(F,E),\Hom_R(E,F)), \quad
T^\pm(X,Y,Z)= XYZ + ZYX.
$$
A Jordan pair that is a sub-pair of such a pair is
called {\it special}.
The Bergman operator is in this case given by
$$
B(X,Y)Z=(\1-XY)Z(\1-YX).
$$
The special case where $V=E \oplus E$ 
gives rise to a self-dual geometry and is
related to the $\K$-Jordan algebra $\End_R(E)$.

\msk
\nin {\bf 8.6. Geometry of right ideals.}
Now let us consider the case of the right $R$-module $V = R$.
In this case (complemented) submodules are the same as
(complemented) right ideals, and the Grassmannian geometry should
be called the {\it geometry of right ideals of $R$}.
Via the bijection $R \to \Hom_R(R,R)$, $r \mapsto l_r$ (left 
multiplication by $r$), the set $\cal P$ of projectors is 
identified with the {\it set of idempotents of $R$},
$$
\Idem(R) := \{ e \in R : \, e^2 = e \}.
$$
The pair $(R,e)$ with an idempotent $e$ is also called a
{\it Morita context} (cf.\ [Lo95, Section 2.1]).
In this case, our Theorem 8.4 corresponds essentially to results
of Loos ([Lo95, Theorem 2.8]).
The symmetric space structure on $\Idem(R)$ 
is described in the same way as after Prop. 8.2:
it is given by
$\mu(e,f) = (2e-\1)f(2e-\1)$.


\msk
\nin {\bf 8.7. Geometry of the projective line.}
Another interesting case is $V=R \oplus R$, taking this decomposition
as base point  $p \in {\cal P}$. 
This gives rise  to a self-dual geometry
(belonging to $R$ seen as a Jordan algebra over $\K$)
which is the {\it projective line over the ring $R$}, see [He95]
and the recent work [BlHa01].  The corresponding $3$-graded
Lie algebra is $\g = \gl_2(R)$, resp., its subalgebra 
$\e_2(R)$ generated by the strict 
upper and lower triangular matrices.

Finally,  let us remark that there exist rings $R$ such that
$R \oplus R \cong R$ as right $R$-modules (e.g., 
take $R=\End_\K(V)$, where $V$ is an infinite dimensional vector
space over a field $\K$; then $V\cong V \oplus V$ as a vector
space, and hence $R=\Hom_\K(V,V) \cong \Hom_\K(V,V \oplus V)$
as a right $R$-module),
 so that the cases 8.6 and 8.7 have non-empty intersection.
\qed

\def\entries{ 

\[Be00 Bertram, W., ``The Geometry of Jordan and Lie Structures,''
Springer LNM {\bf 1754}, Berlin 2000

\[Be01 ---, {\it Generalized projective geometries:
from linear algebra via affine algebra to projective algebra}, 
preprint, Nancy 2001, submitted

\[Be02 ---, {\it Generalized projective geometries: general
theory and equivalence with Jordan structures}, 
Advances in Geometry {\bf 2} (2002), 329--369  

\[Be03 ---, {\it The geometry of null-systems, Jordan algebras
and von Staudt's theorem}, Ann. Inst. Fourier {\bf 53} (1) (2003),
193--225

\[BGN03 Bertram, W., Gl\"ockner, H. and K.-H. Neeb, {\it 
Differential calculus, manifolds and Lie groups over arbitrary
infinite fields}, preprint Darmstadt, Nancy 2003, see
math.GM/0303300

\[BN03 Bertram, W., and K.-H. Neeb, {\it 
Projective completions of Jordan pairs. Part II. 
Manifold structures and symmetric spaces}, in
preparation 

\[BlHa01 Blunck, A. and H. Havlicek, {\it The connected components
of the projective line over a ring}, Advances in Geometry {\bf 1} (2001),
107--117

\[D92 Dorfmeister, J., {\it Algebraic Systems in Differential
Geometry}, p. 9--34 in: ``Jordan Algebras'', Ed. W. Kaup et al.,
de Gruyter, Berlin 1994

\[DNS87 Dorfmeister, J., E. Neher and J. Szmigielski,
{\it Automorphisms of Banach Manifolds associated with the 
$KP$-equation}, Quart. J. Math. Oxford {\bf 40 (2)} (1989), 161--195

\[Fa83 Faulkner, J. R., {\it Stable range and linear groups for
alternative rings}, Geom. Dedicata {\bf 14} (1983), 177--188 


\[Io03  Iordanescu, R., ``Jordan Structures in Geometry and Physics'',
Ed. Ac. Rom\^ane, Bucare\c sti 2003

\[He95 Herzer, A., {\it Chain geometries}. In: F. Buekenhout (editor),
``Handbook of Incidence Geometry'', Elsevier 1995

\[Ka83 Kaup, W., {\it A Riemann Mapping Theorem for Bounded
Symmetric Domains in Complex Banach Spaces},
Math. Z. {\bf 183} (1983), 503--529 

\[Lo67 Loos, O., {\it Spiegelungsr\"aume und homogene symmetrische
R\"aume}, Math. Z. {\bf 99} (1967), 141--170

\[Lo69 ---, ``Symmetric Spaces I,'' Benjamin, New York 1969

\[Lo75 ---, ``Jordan Pairs,'' Springer LNM 460, Berlin 1975

\[Lo77 ---, {\it Bounded Symmetric Domains and Jordan Pairs},
Lecture Notes, Irvine 1977

\[Lo95 ---, {\it Elementary Groups and Stability for Jordan
Pairs}, $K$-Theory {\bf 9} (1995), 77--116

\[Up85 Upmeier, H., ``Symmetric Banach Manifolds and Jordan
$C^*$-algebras,'' North Holland Mathematics Studies, 1985 

}

\references
\dlastpage 

\vfill\eject

\bye 

Alter TExt 8.6, zu (S4):

defines on $\Idem(R)$ a multiplication satisfying the reflection
 space axioms (S1)--(S3). 

Let $\Gr(R)$ denote the set of complemented right ideals of $R$. Then
the inclusion $ {\cal P} \to {\cal C} \times {\cal C} $ from Equation (8.8)
can be seen as an inclusion 
$$
 \Idem(R) \to \Gr(R) \times \Gr(R), \quad 
e \mapsto (eR, (\1-e)R)
 $$
whose image coincides with the set of those pairs 
$(I,J)$ which are transversal in the sense that $I \oplus J = R$.
 In fact, for each such decomposition the 
projection onto $I$ along $J$ is represented by a left multiplication with an idempotent 
$e \in R$. For the tangent spaces we have 
$$ \eqalign{ T_{Re}(\Gr(R)) 
&\cong \Hom_R((\1-e)R, eR) 
\cong \{ \phi \in \Hom_R(R,eR) \cong eR \: \phi(e) = 0\} \cr
&\cong \{ x \in eR \: xe  = 0\} = eR(\1-e) \cr} $$
and likewise 
$$ T_{(\1-e)R}(\Gr(R)) \cong  (\1-e) R e. $$
This leads to 
$$ T_{e}(\Idem(R)) \cong (\1-e)R e \oplus e R(\1-e), $$
and on this space the tangent map (defined in the natural way) 
of $\mu_e$ is $-\id$. In this sense (S4) is also satisfied. 
\qed

\msk
\nin {\bf Further topics.}

-- some more facts on the structure bundle:
{\bf Maybe one should say it already with the definition of $T'X^+$,
that we consider it not as a vector bundle but rather as an affine
bundle: the stabilizer group $P^-$ acts affinely on $V^-$, but this
would not fit with the homogeneous vector bundle picture.!!!} 
note that we can let $P^-$ act in two ways on $T_o'X^+ = \g_{-1}$:
linearly or affinely. The latter acts transitively, hence
$T'X^+$ is homogeneous under $G$. The stabilizer of the base point
``zero vector in $T_o'X^+$'' is then $H$, whence
$T'X^+ \cong G/H =M$ !
Accordingly, we have an embedding of $T'X^+$ into $X^+ \times X^-$,
given by ...
(in the real non-deg. case: $M$ is isomorphic to the
``cotangent bundle of $X^+$".)
Next: definition of generalized cross-ratios and Maslov-indices
(the latter seem to need the existence of a polarity).
Natural: to a quadruple $(x,a,y,b)$ in generic position
associate the operator
$$
K_{x,a} K_{y,a}^{-1} K_{y,b} K_{x,b}^{-1}:T_x X^+ \to T_x X^+.
$$
(write as diagram, closed rectangle; of course we could start at
any corner of this rectangle).
This defines a $G$-invariant section over generic points in
$X^+ \times X^- \times X^+ \times X^-$. It seems natural to
relate it to some fundamental group, seeing the rectangle
as a loop attached to $x$. 

-- it seems tempting to try a similar theory for $5$-gradings.
In particular, if $(e^+,e^-)$ is a pair of idempotents {\rm(of what???)}, the
derivation $T(e^+,e^-)$ of $\g$ does induce a $5$-grading
(but so far I do not have a good conceptual proof of this).
This should be a good point to organize the Pierce-theory.

{\bf Associative algebras:}

\Lemma 7.4. If $A$ is an associative unital algebra and $a, b \in A$ commute, then 
$ab \in A^\times$ is equivalent to $a,b \in A^\times$. 

\Proof. If $a$ and $b$ are invertible, then $ab$ is invertible with inverse 
$b^{-1} a^{-1} = a^{-1} b^{-1}$. If, conversely, $ab$ is invertible and 
$c$ is an inverse of $ab$, then $a$ and $b$ commute with $c$. Hence 
$c a$ is an inverse of $b$ and $cb$ is an inverse of $a$. 
\qed

\Examples 7.5. (a) If $V$ and $W$ are non-zero vector spaces and 
$A := \End(V \oplus W)$, then the projection $e$ of $V \oplus W$ onto $V$ is an idempotent, 
and the corresponding Morita context leads to the Jordan pair 
$(V^+, V^-)$ given by 
$$ V^+ := \Hom(W,V) \quad \hbox{ and } \quad 
 V^- := \Hom(V,W). $$

For $x \in \End(W)$ let $L(x)$ denote the left multiplication with $x$ on $\Hom(V,W)$ 
and for $y \in \End(V)$ let $R(y)$ denote the right multiplication with $y$. 
Then the Bergman operator of the Jordan pair $(V^+, V^-)$ is given by 
$$ B(x,y) = L(\1 - xy) R(\1 - yx). $$
Therefore Lemma 7.4 implies that 
the pair $(x,y)$ is quasi-invertible if and only if the operators 
$L(\1 - xy)$ and $R(\1 - yx)$ are invertible. This in turn is equivalent to 
$\1 - xy \in \GL(W)$ and $\1 - yx \in \GL(V)$. 
If $\1 -xy \in \GL(W)$ and $z := (\1 - xy)^{-1} \in \GL(W)$, then 
$\1 + yzx$ is an inverse of $\1 - yx$ in $\End(V)$. Therefore the invertibility of 
$\1 - xy$ implies the invertibility of $\1 - yx$ and vice versa. Therefore 
$$ (x,y) \hbox{ \ quasiinvertible} \quad 
\Longleftrightarrow \quad \1 - xy \in \GL(W)\quad 
\Longleftrightarrow \quad \1 - yx \in \GL(V). $$

\par\nin (b) A similar argument as in (a) implies that for the Jordan pair 
$(A,A)$ associated to an associative algebra $A$, the pair $(x,y)$ is quasi-invertible 
if and only if $\1 - xy \in A^\times$ if and only if $\1 - yx \in
A^\times$. 
This is a special case of the {\it symmetry principle} for Jordan
pairs (cf.\ [Lo95]). 
\qed

\bye